# ASYMPTOTICS FOR POSTERIOR HAZARDS

By Pierpaolo De Blasi,[1] Giovanni Peccati[2] and Igor Prünster[1]

*Università di Torino, Collegio Carlo Alberto, Université Paris Ouest and Université Paris VI and Università di Torino and Collegio Carlo Alberto and ICER*

An important issue in survival analysis is the investigation and the modeling of hazard rates. Within a Bayesian nonparametric framework, a natural and popular approach is to model hazard rates as kernel mixtures with respect to a completely random measure. In this paper we provide a comprehensive analysis of the asymptotic behavior of such models. We investigate consistency of the posterior distribution and derive fixed sample size central limit theorems for both linear and quadratic functionals of the posterior hazard rate. The general results are then specialized to various specific kernels and mixing measures yielding consistency under minimal conditions and neat central limit theorems for the distribution of functionals.

**1. Introduction.** Bayesian nonparametric methods have found a fertile ground of applications within survival analysis. Indeed, given that survival analysis typically requires function estimation, the Bayesian nonparametric paradigm seems to be tailor made for such problems, as already shown in the seminal papers by Doksum [4], Dykstra and Laud [6], Lo and Weng [24] and Hjort [11]. According to the approach of [6, 24], the hazard rate is modeled as a mixture of a suitable kernel with respect to an increasing additive process (see [32]) or, more generally, a completely random measure (see [21]). This approach will be the focus of the present paper: below we first present the model and, then, the two asymptotic issues we are going to tackle, namely weak consistency and the derivation of fixed sample size central limit theorems (CLTs) for functionals of the posterior hazard rate.

Received February 2008; revised June 2008.

[1]Supported in part by MIUR, Grant 2006/133449.

[2]Supported in part by ISI Foundation, Lagrange Project.

1.1. *Life-testing model with mixture hazard rate.* Denote by $Y$ a positive absolutely continuous random variable representing the lifetime and assume that its random hazard rate is of the form

$$(1) \qquad \tilde{h}(t) = \int_{\mathbb{X}} k(t, x) \tilde{\mu}(dx),$$

where $k$ is a kernel and $\tilde{\mu}$ a completely random measure on some Polish space $\mathbb{X}$ endowed with its Borel $\sigma$-field $\mathscr{X}$. The kernel $k$ is a jointly measurable application from $\mathbb{R}^+ \times \mathbb{X}$ to $\mathbb{R}^+$ and the application $C \mapsto \int_C k(t, x) \, dt$ defines a $\sigma$-finite measure on $\mathscr{B}(\mathbb{R}^+)$ for any $x$ in $\mathbb{X}$. Typical choices, which we will also consider in this paper, are:

(i) the Dykstra–Laud (DL) kernel [6]

$$(2) \qquad k(t, x) = \mathbb{I}_{(0 \le x \le t)},$$

which leads to monotone increasing hazard rates;

(ii) the rectangular kernel (see, e.g., [13]) with bandwidth $\tau > 0$

$$(3) \qquad k(t, x) = \mathbb{I}_{(|t-x| \le \tau)};$$

(iii) the Ornstein–Uhlenbeck (OU) kernel (see, e.g., [25, 26]) with $\kappa > 0$

$$(4) \qquad k(t, x) = \sqrt{2\kappa} \exp(-\kappa(t-x)) \mathbb{I}_{(0 \le x \le t)};$$

(iv) the exponential kernel (see, e.g., [14])

$$(5) \qquad k(t, x) = x^{-1} e^{-t/x},$$

which yields monotone decreasing hazard rates.

As for the mixing measure in (1), letting $(\mathbb{M}, \mathscr{B}(\mathbb{M}))$ be the space of boundedly finite measures on $(\mathbb{X}, \mathscr{X})$, $\tilde{\mu}$ is taken to be a *completely random measure* (CRM) in the sense of [21]. This means that $\tilde{\mu}$ is a random element defined on $(\Omega, \mathscr{F}, \mathbb{P})$, taking values in $(\mathbb{M}, \mathscr{B}(\mathbb{M}))$ and such that, for any collection of disjoint sets, $B_1, B_2, \ldots$, the random variables $\tilde{\mu}(B_1), \tilde{\mu}(B_2), \ldots$ are mutually independent. Appendix A.1 provides a brief account of CRMs, as well as justifications of the following statements. It is important to recall that a CRM is characterized by its Poisson intensity $\nu$, which we can write as

$$(6) \qquad \nu(dv, dx) = \rho(dv | x) \lambda(dx),$$

where $\lambda$ is a $\sigma$-finite measure on $\mathbb{X}$. If, furthermore, $\nu(dv, dx) = \rho(dv)\lambda(dx)$, the corresponding CRM $\tilde{\mu}$ is termed *homogeneous*, otherwise it is said to be *nonhomogeneous*. We always consider kernels such that $\int_{\mathbb{X}} k(t, x)\lambda(dx) < +\infty$. Throughout the paper, we will take $\nu$ and $\lambda$ to be nonatomic and we shall moreover assume that

$$(\text{H1}) \qquad \rho(\mathbb{R}^+ | x) = +\infty \qquad \text{a.e.-}\lambda \quad \text{and} \quad \text{supp}(\lambda) = \mathbb{X},$$

where supp($\tau$) indicates the topological support of a given measure $\tau$. Note that (H1) is equivalent to requiring that $\tilde{\mu}$ jumps infinitely often on any bounded set of positive $\lambda$-measure and is indeed a desirable property for a mixing measure, since it ensures that the topological support of $\tilde{\mu}$ is the whole space $\mathbb{M}$. See also the discussion around formula (3.22) in [18] for an account of the usefulness of (H1) for inferential purposes. In the examples we will focus on a large class of CRMs, which includes almost all CRMs used so far in applications and is characterized by an intensity measure of the type

$$(7) \qquad \nu(dv, dx) = \frac{1}{\Gamma(1-\sigma)} \frac{e^{-\gamma(x)v}}{v^{1+\sigma}} \, dv \lambda(dx),$$

where $\sigma \in [0, 1)$ and $\gamma$ is a strictly positive function on $\mathbb{X}$. Note that, if $\gamma$ is a constant, the resulting CRMs coincide with the generalized gamma measures [2], whereas when $\sigma = 0$ they are extended gamma CRMs [6, 24].

Having defined the ingredients of the mixture hazard (1), we can complete the description of the model, which is often referred to as *life-testing model*. The cumulative hazard is then given by $\tilde{H}(t) = \int_0^t \tilde{h}(s) \, ds$ and, provided

$$(8) \qquad \tilde{H}(t) \to \infty \qquad \text{for } t \to +\infty \qquad \text{a.s.},$$

one can define a random density function $\tilde{f}$ as

$$(9) \qquad \tilde{f}(t) = \tilde{h}(t) \exp(-\tilde{H}(t)) = \tilde{h}(t)\tilde{S}(t),$$

where $\tilde{S}(t) := \exp(-\tilde{H}(t))$ is the *survival function*, providing the probability that $Y > t$. Consequently, the random cumulative distribution function of $Y$ is of the form $\tilde{F}(t) = 1 - \exp(-\tilde{H}(t))$. Note that, given $\tilde{\mu}$, $\tilde{h}$ represents the hazard rate of $Y$, that is, $h(t) \, dt = \mathbb{P}(t \leq Y \leq t + dt | Y \geq t, \tilde{\mu})$. Throughout the paper we will assume that

$$(H2) \qquad \mathbb{E}[\tilde{H}(t)] = \int_0^t \int_{\mathbb{R}^+ \times \mathbb{X}} v k(u, x) \rho(dv|x) \lambda(dx) \, du < +\infty \qquad \forall t > 0.$$

Such models have recently received much attention due to their relatively simple implementation in applications. Important developments, dealing also with more general multiplicative intensity models, can be found in [12, 13, 14, 15, 25, 26], among others.

1.2. *Posterior consistency.* The study of consistency of Bayesian nonparametric procedures represents one of the main recent research topics in Bayesian theory. The "frequentist" (or "what if") approach to Bayesian consistency consists of generating independent data from a "true" fixed density $f_0$ and checking whether the sequence of posterior distributions accumulates in some suitable neighborhood of $f_0$. Specifically, denote by $P_0$ the probability distribution associated with $f_0$ and by $P_0^\infty$ the infinite product measure.

Moreover, the symbol $\mathbb{F}$ indicates the space of density functions absolutely continuous with respect to the Lebesgue measure on $\mathbb{R}$, endowed with the Borel $\sigma$-field $\mathscr{B}(\mathbb{F})$ (with respect to an appropriate $L^1$-topology). Now, if $\Pi$ is the prior distribution of some random density function $\tilde{f}$, taking values in $\mathbb{F}$, and $\Pi_n$ denotes its posterior distribution, then one is interested in establishing sufficient conditions to have that, as $n \to +\infty$, for any $\varepsilon > 0$

$$\text{(10)} \qquad \Pi_n(A_\varepsilon(f_0)) \to 1 \qquad \text{a.s.-}P_0^\infty,$$

where $A_\varepsilon(f_0)$ represents a $\varepsilon$-neighborhood of $f_0$ in a suitable topology. If (10) holds, then $\Pi$ is said to be *consistent* at $f_0$. Now, if $A_\varepsilon(f_0)$ is chosen to be a weak neighborhood, one obtains *weak consistency*. Sufficient conditions for weak consistency of various important nonparametric models have been provided in, for example, [8, 33, 35, 37]. By requiring (10) to hold with $A_\varepsilon$ being a $L_1$-neighborhood, one obtains the stronger notion of $L_1$ *consistency*: general sufficient conditions for this to happen are provided in [1, 8, 36]. In the context of discrete models such as neutral to the right processes, posterior consistency has been studied in [9, 19, 20]. For a thorough review of the literature on consistency issues, the reader is referred to the monograph [10].

Turning back to the life-testing model defined by (1) and (9), little is known about consistency, since their structure is intrinsically very different from the models considered so far. First results were given in [5, 25]. In particular, in [5] consistency is established for the DL kernel with extended gamma mixing measure assuming a bounded "true" hazard. In this paper, we determine sufficient conditions for weak consistency of Bayesian nonparametric models defined in terms of mixture random hazard rates. We also cover the case of lifetimes subject to independent right-censoring. Then, we use this general result for establishing weak consistency for mixture hazards with the specific kernels in (2)–(5) and CRMs characterized by (7). In particular, we obtain consistency essentially w.r.t. nondecreasing hazards for DL mixtures, w.r.t. bounded Lipschitz hazards for rectangular mixtures, w.r.t. to hazards with certain local exponential decay rate for OU mixtures and w.r.t. completely monotone hazards for exponential mixtures.

1.3. *Functionals of the posterior mixture hazard rate.* The second aspect we investigate is the asymptotic behavior (in the sense of larger and larger time horizons) of functionals of the posterior random hazard rate given a fixed number of observations. We shall focus on functionals of statistical relevance, such as means, path-second moments and path-variances. Indeed, any CLT involving this type of functionals may be used to derive a synthetic—yet highly informative—picture of the "global shape" of a given (prior or posterior) hazard rate model. In particular, as we will see below, CLTs for linear and quadratic functionals contain specific information about the trend, the oscillations and the overall asymptotic variance of a random object such as (1). This represents an important issue since, though widely used in practice, the implications of the choice of specific kernels and CRMs in defining (1) are generally not well understood and their choice is based on mere empirical considerations. In [27] functionals of the prior hazard rate are considered: the results, despite being of theoretical relevance, can serve also as a guide for prior specification. For instance, it is shown that the trend of the cumulative hazard with a DL kernel (2) with a homogeneous CRM is $T^2$, with the oscillations around the trend increasing like $T^{3/2}$, whereas with a rectangular kernel the trend is $T$ and the oscillations increase like $T^{1/2}$. Moreover, the parameters of the kernel and the CRM enter the variance of the asymptotic Gaussian random variable, thus leading to a rigorous procedure for their a priori selection.

Here, we face the more challenging problem of deriving CLTs for the posterior hazard rate: indeed, the model defined by (1) and (9) is not conjugate and, hence, the derivation of distributional results for posterior functionals is quite demanding. However, by exploiting the posterior representation of James [15] (to be detailed in Section 2), we are able to provide fixed sample size CLTs also for functionals of posterior hazard rates. One of our main findings is that, in all the considered special cases, the CLTs associated with the posterior hazard rate are *the same as for the prior ones*, and this for any number of observations. If one interprets CLTs as approximate "global pictures" of a model, the conclusions to be drawn from our results are quite clear. Indeed, although consistency implies that a given model can be asymptotically directed toward any deterministic target, the overall structure of a posterior hazard rate is systematically determined by the prior choice, even after conditioning on a very large number of observations.

As an example of the results derived in the sequel, consider again the hazard rate given by the DL kernel (2) with a homogeneous CRM, and let $\mathbf{Y} = (Y_1, \ldots, Y_n)$ be a set of observations. In Section 4.3.1, we will prove that

$$T^{-3/2}[\tilde{H}(T) - cT^2] | \mathbf{Y} \overset{\text{law}}{\longrightarrow} X$$

(the precise meaning of such a conditional convergence in law will be clarified in the sequel), where $c$ is a constant and $X$ is a centered Gaussian random variable with variance $\sigma^2$. As anticipated, the crucial point will be that both $c$ and $\sigma^2$ are independent of $n$ and $\mathbf{Y}$, and that they are actually the same constants appearing in the prior CLTs proved in [27]. A more detailed illustration of these phenomena is provided in Section 4.3, where we also discuss analogous results involving other models, as well as limit theorems for quadratic functionals.

We stress that our choice of $+\infty$ as a limiting point is mainly conventional, and that one can easily modify our framework to deal with models that live within a finite window of time by using an appropriate deformation of the time scale. For instance, one can embed a hazard rate model defined on $[0, +\infty)$ into a finite time interval, by substituting the time parameter $T$ in the previous discussion with an increasing function of the type $\log[T^*/(T^* - T)]$, where $T^* < +\infty$ and $0 \leq T < T^*$.

1.4. *Outline.* The paper is organized as follows. Section 2 provides the posterior characterization of model (1). In Section 3 sufficient conditions for weak consistency are established. Section 4 deals with posterior linear and quadratic functionals of the mixture hazard. The results are illustrated by various examples involving specific kernels and CRMs. In Section 5 some concluding remarks and future research lines are presented. Further results, which are also of independent interest, and the proofs are deferred to the Appendix.

## 2. Posterior distribution of the random hazard rate.

In order to make Bayesian inference starting from model (1), an explicit posterior characterization is essential. Indeed, the first treatments of model (1) were limited to considering extended gamma CRMs, which allow for a relatively simple posterior characterization [6, 24]. Analysis beyond gamma-like choices of $\tilde{\mu}$ has not been possible for a long time due to the lack of a suitable and implementable posterior characterization: however, in James [15] this goal has been achieved and many choices for $\tilde{\mu}$ can now be explored. See also [23] for a different derivation of these results. In what follows, we give an explicit description of the posterior characterization of the model (1).

Let $\tilde{P}_{\tilde{f}}$ be the random probability measure associated with (9) and denote by $(Y_n)_{n \geq 1}$ a sequence of exchangeable observations, defined on $(\Omega, \mathscr{F}, \mathbb{P})$ and taking values in $\mathbb{R}^+$, such that, given $\tilde{P}_{\tilde{f}}$, the $Y_n$'s are i.i.d. with distribution $\tilde{P}_{\tilde{f}}$, that is, $\mathbb{P}[Y_1 \in B_1, \ldots, Y_n \in B_n | \tilde{P}_{\tilde{f}}] = \prod_{i=1}^n \tilde{P}_{\tilde{f}}(B_i)$ for any $B_i \in \mathscr{B}(\mathbb{R}^+)$, $i = 1, \ldots, n$ and $n \geq 1$. The joint (conditional) density of $\mathbf{Y} = (Y_1, \ldots, Y_n)$ given $\tilde{\mu} = \mu$ is then given by

$$e^{-\int_{\mathbb{X}} \sum_{i=1}^n \int_0^{y_i} k(t,x)\,dt\,\mu(dx)} \prod_{i=1}^n \int_{\mathbb{X}} k(y_i, x)\mu(dx).$$

In this context it is important to consider also some censoring mechanism, specifically independent right-censoring. Hence, suppose there are additionally $Y_{n+1}, \ldots, Y_m$ random times which are right censored by censoring times $C_{n+1}, \ldots, C_m$, that is, $Y_i > C_i$ for $i = n+1, \ldots, m$ [by exchangeability, it would be equivalent to assume the right censored data to be an arbitrary $(m-n)$-dimensional subvector of $(Y_1, \ldots, Y_m)$]. It is well known that assuming the distribution of $C$ to be known is equivalent to assuming the distribution of $C$ is a priori independent of the distribution of $Y$. Hence, the posterior distribution of $\tilde{\mu}$ may be obtained without even specifying the prior on the distribution of $C$. Then the likelihood function based on $\mathbf{Y} = (Y_1, \ldots, Y_m)$, where the vector $\mathbf{Y}$ is composed of $n$ completely observed times and $m - n$ right censored times, has the form

$$(11) \qquad \mathscr{L}(\mu; \mathbf{y}) = e^{-\int_{\mathbb{X}} K_m(x)\mu(dx)} \prod_{i=1}^{n} \int_{\mathbb{X}} k(y_i, x)\mu(dx),$$

where $K_m(x) = \sum_{i=1}^{m} \int_{0}^{y_i \wedge c_i} k(t, x)\, dt$ and we set $c_i = \infty$ for $i = 1, \ldots, n$. If we now augment the likelihood with respect to the latent variables $\mathbf{X} = (X_1, \ldots, X_n)$, (11) reduces to

$$\mathscr{L}(\mu; \mathbf{y}, \mathbf{x}) = e^{-\int_{\mathbb{X}} K_m(x)\mu(dx)} \prod_{i=1}^{n} k(y_i; x_i)\mu(dx_i)$$

$$= e^{-\int_{\mathbb{X}} K_m(x)\mu(dx)} \prod_{j=1}^{k} \mu(dx_j^*)^{n_j} \prod_{i \in D_j} k(y_i; x_j^*),$$

where $\mathbf{X}^* = (X_1^*, \ldots, X_k^*)$ denote the $k \leq n$ distinct latent variables, $n_j$ is the frequency of $X_j^*$ and $D_j = \{r : x_r = x_j^*\}$. Finally, set $\tau_{n_j}(x) = \int_{\mathbb{R}^+} v^{n_j} e^{-vK_m(x)} \rho(dv|x)$. We are now in a position to state the posterior characterization of the mixture hazard rate.

THEOREM 1 (James [15]). *Let $\tilde{h}$ be a random hazard rate as defined in (1), corresponding to model (9). Then, given $\mathbf{Y}$, the posterior distribution of $\tilde{h}$ can be characterized as follows:*

(i) *Given $\mathbf{X}$ and $\mathbf{Y}$, the conditional distribution of $\tilde{\mu}$ coincides with the distribution of the random measure*

$$(12) \qquad \tilde{\mu}^{m,*} + \sum_{i=1}^{k} J_i \delta_{X_j^*} = \tilde{\mu}^{m,*} + \Delta^{n,*},$$

*where $\tilde{\mu}^{m,*}$ is a CRM with intensity measure*

$$(13) \qquad \nu^{m,*}(dv, dx) := e^{-vK_m(x)}\rho(dv|x)\lambda(dx),$$

$\Delta^{n,*}(dx) := \sum_{i=1}^{k} J_i \delta_{X_i^*}(dx)$ *with, for $i = 1, \ldots, k$, $X_i^*$ a fixed point of discontinuity with corresponding jump $J_i$ distributed as*

$$(14) \qquad f_{J_i}(dv) = \frac{v^{n_i} e^{-vK_m(X_i^*)} \rho(dv|X_i^*)}{\int_{\mathbb{R}^+} v^{n_i} e^{-vK_m(X_i^*)} \rho(dv|X_i^*)}.$$

*Moreover, the $J_i$'s are, conditionally on $\mathbf{X}$ and $\mathbf{Y}$, independent of $\tilde{\mu}^{m,*}$.*

(ii) *Conditionally on* $\mathbf{Y}$, *the distribution of the latent variables* $\mathbf{X}$ *is*

$$f(dx_1^*, \ldots, dx_k^* | \mathbf{Y}) = \frac{\prod_{j=1}^k \tau_{n_j}(x_j^*) \prod_{i \in D_j} k(y_i, x_j^*) \lambda(dx_j^*)}{\sum_{k=1}^n \sum_{\mathbf{n} \in A_{k,n}} \prod_{j=1}^k \int_{\mathbb{X}} \tau_{n_j}(x) \prod_{i \in D_j} k(y_i, x) \lambda(dx)}$$

*for any* $k \in \{1, \ldots, n\}$ *and* $\mathbf{n} := (n_1, \ldots, n_k) \in A_{k,n} := \{(n_1, \ldots, n_k) : n_j \geq 1, \sum_{j=1}^k n_j = n\}$.

**3. Consistency.** Our first goal consists in deriving sufficient conditions for weak consistency of the Bayesian nonparametric life-testing model (9) with mixture hazard (1), which covers also the case of data subject to right-censoring. Then, we exploit this criterion for obtaining consistency results for specific mixture hazards.

In the case of complete data, a general and widely used sufficient condition for *weak consistency* with respect to a "true" unknown density function $f_0$, due to Schwartz [33], requires a prior $\Pi$ to assign positive probability to Kullback–Leibler neighborhoods of $f_0$, that is,

(15)           $\Pi(f \in \mathbb{F} : d_{KL}(f_0, f) < \varepsilon) > 0$    for any $\varepsilon > 0$,

where $d_{KL}(f_0, f) = \int \log(f_0(t)/f(t)) f_0(t) \, dt$ denotes the Kullback–Leibler divergence between $f_0$ and $f$.

In the presence of right-censoring, we do not actually observe the lifetime $Y$, but, $(Z, \Delta)$, where $Z = Y \wedge C$, $\Delta = \mathbb{I}_{(Y \leq C)}$ for $C$ a censoring time with distribution $P_c$ admitting density $f_c$. Clearly, this leads us to consider a prior on the space $\mathbb{F} \times \mathbb{F}$ and the corresponding prior $\Pi^*$ induced on the space of the distribution of the observables $(Z_i, \Delta_i)$'s.

The strategy of the proof consists in first rewriting the Kullback–Leibler condition in terms of the induced prior $\Pi^*$: this condition then guarantees consistency of $\Pi^*$. Moreover, it allows us to deduce the consistency of $\Pi$, the prior on the distribution of the lifetime $Y$, under independent right-censoring with the simple support condition

(16)                          $\mathrm{supp}(P_c) = \mathbb{R}^+.$

The last step consists in translating the Kullback–Leibler condition into a condition in terms of uniform neighborhoods of the true hazard rate $h_0$ on the interval $(0, T]$ for any finite $T$. When dealing with models for hazard rates, the latter appears to be both more natural and easy to verify.

Without risk of confusion, in the following we denote by $\Pi$ the prior on $\tilde{f}$ and also the prior induced on $\tilde{h}$. Moreover, recall that the "true" density $f_0$ can always be represented in terms of the "true" hazard $h_0$ as $f_0(t) = h_0(t) \exp(-\int_0^t h_0(s) \, ds)$.

THEOREM 2. *Let $\tilde{f}$ be a random density function defined by (1) and (9) with kernels (2)–(5) and denote its (prior) distribution by $\Pi$. Suppose the distribution of the censoring times $P_c$ is independent of the lifetime $Y$, absolutely continuous and satisfies (16). Moreover, assume that the following conditions hold:*

(i) *$f_0(t)$ is strictly positive on $(0, \infty)$ and $\int_{\mathbb{R}^+} \max\{\mathbb{E}[\tilde{H}(t)], t\} f_0(t)\, dt < \infty$;*

(ii) *there exists $r > 0$ such that $\liminf_{t \downarrow 0} \tilde{h}(t)/t^r = \infty$ a.s.*

*Then, a sufficient condition for $\Pi$ to be weakly consistent at $f_0$ is that*

$$(17) \qquad \Pi\left\{ h : \sup_{0 < t \leq T} |h(t) - h_0(t)| < \delta \right\} > 0$$

*for any finite $T$ and positive $\delta$.*

Some comments regarding the conditions are in order at this point. Let us start by condition (i): the strict positivity of $f_0$ on $(0, \infty)$ is equivalent to strict positivity of the "true" hazard $h_0$ on $(0, \infty)$, which is a property satisfied by any reasonable $h_0$. The second part of condition (i), which is also related to the asymptotic characterizations considered in Section 4, clearly becomes more restrictive the faster the trend of the cumulative hazard. However, note that if $h_0$ is a power function, then $f_0$ admits moments of any order and, hence, it is enough that the trend of the cumulative hazard is a power function as well. Condition (ii) allows to remove the somehow artificial assumption of $h_0(0) > 0$ as in [5]. Indeed, $h_0(0) = 0$ represents a common situation in practice and condition (ii) covers such a case by controlling the small time behavior of $\tilde{h}$. Obviously, if $h_0(0) > 0$, then one would adopt a random hazard $\tilde{h}$ nonvanishing in 0 and so condition (ii) would be automatically satisfied. Overall, the result can be seen as a general consistency criterion for mixture hazard models and deals automatically with the case of independent right-censoring. Moreover, it should be extendable in a quite straightforward way to mixture hazards with different reasonably behaving kernels.

Before entering a detailed analysis of specific models, we show how condition (ii) of Theorem 2 can be reduced to the problem of studying the short time behavior of the CRM and, moreover, we establish that the CRMs defined in (7) satisfy the corresponding short time behavior requirement. Throughout this section we assume $\mathbb{X} = \mathbb{R}^+$ and, hence, when useful, $\tilde{\mu}$ will be treated as an increasing additive process (see [32]), namely the càdlàg distribution function induced by $\tilde{\mu}$.

PROPOSITION 3. *Let $\tilde{h}$ be a mixture hazard (1). Then condition (ii) in Theorem 2 is implied by:*

(ii1) *there exists* $\varepsilon > 0$ *such that* $\tilde{h}(t) \geq c\tilde{\mu}((0,t])$ *for* $t < \varepsilon$, *where* $c$ *is a constant not depending on* $t$;

(ii2) *there exists* $r > 0$ *such that* $\liminf_{t\downarrow 0} \tilde{\mu}((0,t])/t^r = \infty$ *a.s.*

*In particular,* (ii1) *holds if* $k$ *is either the DL* (2) *or the OU* (4) *kernel;* (ii2) *holds if* $\tilde{\mu}$ *is a CRM belonging to* (7) *with* $\sigma \in (0,1)$ *and* $\lambda(dx) = dx$.

Condition (ii1) requires that the random hazard leaves the origin at least as fast as the driving CRM, which is typically the case. Out of the four considered kernels, we have to face the problem of $\tilde{h}(0) = 0$ a.s. for the DL and OU mixtures and for both kernels (ii1) is satisfied. Condition (ii2) asks to control the small time behavior of the CRM and is met by CRMs like (7). If one is interested in CRMs different from (7), one can try to adapt one of the several results on small time behavior known in the literature (see, e.g., [32] and references therein).

We now move on to deriving explicit consistency results for mixture hazard life-testing models based on the four kernels defined in (2)–(5). These results are derived by verifying the conditions of Theorem 2 and, thus, hold also for data subject to right-censoring with absolutely continuous censoring distribution satisfying (16). Though the details of the proofs are different, they rely on a common strategy: first consistency is established via condition (17) for "true" hazards of mixture form $h_0(t) = \int_{\mathbb{R}^+} k(t,x)\mu_0(dx)$, where $k$ is the same kernel used for defining the specific model $\tilde{h}$; then, we show that these mixture $h_0$'s are arbitrarily close in the uniform metric to any $h_0$ belonging to a class of hazards having a suitable qualitative feature.

We first deal with DL mixture hazards $\tilde{h}(t) = \int_{\mathbb{R}^+} \mathbb{I}_{(0 \leq x \leq t)} \tilde{\mu}(dx)$, which represent a model for nondecreasing hazard rates. The result establishes weak consistency of such models for any nondecreasing $h_0$ satisfying some mild additional conditions.

THEOREM 4. *Let* $\tilde{h}$ *be a mixture hazard* (1) *with DL kernel and* $\tilde{\mu}$ *satisfying condition* (ii2) *of Proposition 3.*

*Then* $\Pi$ *is weakly consistent at any* $f_0 \in \mathscr{F}_1$, *where* $\mathscr{F}_1$ *is defined as the set of densities for which:* (i) $\int_{\mathbb{R}^+} \mathbb{E}[\tilde{H}(t)]f_0(t)\,dt < \infty$; (ii) $h_0(0) = 0$ *and* $h_0(t)$ *is strictly positive and nondecreasing for any* $t > 0$.

The second model we consider is represented by rectangular mixture hazards $\tilde{h}(t) = \int_{\mathbb{R}^+} \mathbb{I}_{(|t-x| \leq \tilde{\tau})} \tilde{\mu}(dx)$. In order to obtain consistency with respect to a large class of $h_0$'s we treat the bandwidth $\tau$ as a hyper-parameter and assign to it an independent prior $\pi$, whose support contains $[0, L]$ for some $L > 0$. So we have two sources of randomness: $\tilde{\tau}$ with distribution $\pi$ and $\tilde{\mu}$, whose distribution we denote by $Q$. Hence, the prior distribution $\Pi$ on $\tilde{h}$ is induced by $\pi \times Q$ via the map $(\tau, \mu) \to h(\cdot|\tau, \mu) := \int \mathbb{I}_{(|\cdot-x| \leq \tau)} \mu(dx)$. In this framework we are able to derive consistency at essentially any bounded and nonvanishing Lipschitz hazard $h_0$.

THEOREM 5. *Let $\tilde{h}$ be a mixture hazard (1) with rectangular kernel and random bandwidth $\tilde{\tau}$ independent of $\tilde{\mu}$. Moreover, the support of the prior $\pi$ on $\tilde{\tau}$ contains $[0, L]$ for some $L > 0$.*

*Then $\Pi$ is weakly consistent at any $f_0 \in \mathscr{F}_2$, where $\mathscr{F}_2$ is defined as the set of densities for which:* (i) $\int_{\mathbb{R}^+} \max\{\mathbb{E}[\tilde{H}(t)], t\} f_0(t) \, dt < \infty$; (ii) $h_0(t) > 0$ *for any $t \geq 0$;* (iii) $h_0$ *is bounded and Lipschitz.*

Now consider OU mixture hazards $\tilde{h}(t) = \int_{\mathbb{R}^+} \sqrt{2\kappa} e^{-\kappa(t-x)} \mathbb{I}_{(0 \leq x \leq t)} \tilde{\mu}(dx)$. Define for any differentiable decreasing function $g$ the local exponential decay rate as $-g'(y)/g(y)$. Our result establishes consistency at essentially any $h_0$ which exhibits, in regions where it is decreasing, a local exponential decay rate smaller than $\kappa\sqrt{2\kappa}$. This sheds also some light on the role of the kernel-parameter $\kappa$: choosing a large $\kappa$ leads to less smooth trajectories of $\tilde{h}$, but, on the other hand, ensures also consistency with respect to $h_0$'s which have abrupt decays in certain regions.

THEOREM 6. *Let $\tilde{h}$ be a mixture hazard (1) with OU kernel and $\tilde{\mu}$ satisfying condition* (ii2) *of Proposition 3.*

*Then $\Pi$ is weakly consistent at any $f_0 \in \mathscr{F}_3$, where $\mathscr{F}_3$ is defined as the set of densities for which:* (i) $\int_{\mathbb{R}^+} \max\{\mathbb{E}[\tilde{H}(t)], t\} f_0(t) \, dt < \infty$; (ii) $h_0(0) = 0$ *and $h_0(t) > 0$ for any $t > 0$;* (iii) $h_0$ *is differentiable and, for any $t > 0$ such that $h'_0(t) < 0$, the corresponding local exponential decay rate is smaller than $\kappa\sqrt{2\kappa}$.*

REMARK 1. In the above three mixture hazard models, one typically selects CRMs with $\lambda$ in (6) being the Lebesgue measure on $\mathbb{R}^+$. If this is the case, then condition (i) in the definition of $\mathscr{F}_i$ $(i = 1, 2, 3)$, becomes $\int_{\mathbb{R}^+} t^2 f_0(t) \, dt < \infty$ for DL mixture hazards and $\int_{\mathbb{R}^+} t f_0(t) \, dt < \infty$ for rectangular and OU mixtures.

Now we deal with mixture hazards based on an exponential kernel $\tilde{h}(t) = \int_{\mathbb{R}^+} x^{-1} e^{-t/x} \tilde{\mu}(dx)$, which are used to model decreasing hazard rates. Note that, in contrast to the DL, rectangular and OU kernels which all exhibit, for any fixed $t$, finite support on $\mathbb{R}^+$ when seen as functions of $x$, in this case the support is $\mathbb{R}^+$ for any fixed $t$. This implies the need for quite different techniques for handling it. Recall that a function $g$ on $\mathbb{R}^+$ is completely monotone if it possesses derivatives $g^{(n)}$ of all orders and $(-1)^n g^{(n)}(y) \geq 0$ for any $y > 0$. The next result shows that consistency holds at essentially any completely monotone hazard for which $h_0(0) < \infty$.

THEOREM 7. *Let $\tilde{h}$ be a mixture hazard (1) with exponential kernel such that $\tilde{h}(0) < \infty$ a.s.*

*Then $\Pi$ is weakly consistent at any $f_0 \in \mathscr{F}_4$, where $\mathscr{F}_4$ is defined as the set of densities for which:* (i) $\int_{\mathbb{R}^+} t f_0(t) \, dt < \infty$; (ii) $h_0(0) < \infty$; (iii) $h_0$ *is completely monotone.*

Note that the requirement of $\tilde{h}$ not to explode in 0 is easily achieved by selecting $\lambda$ in (6) such that $\int_{\mathbb{R}^+ \times \mathbb{R}^+} (1 - e^{-ux^{-1}v}) \rho(dv|x) \lambda(dx) < \infty$ for all $u > 0$, which is equivalent to $\tilde{h}(0) < \infty$ a.s. [see (36) in Appendix A.1].

## 4. Fixed sample size posterior CLTs.
In this section we derive CLTs for functionals of the random hazard given a fixed set of observations as time diverges. For the sake of clarity, in the following we confine ourselves to the case of complete observations; however, all subsequent results immediately carry over to the case of data subject to right-censoring.

### 4.1. *Further concepts and notation.*
Since we will heavily exploit the posterior characterization of $\tilde{h}$ recalled in Theorem 1, it is useful to introduce first some definitions related to quantities involved in its statement. Whenever convenient, we shall use the notation $\nu^{0,*} := \nu$ and $\tilde{\mu}^{0,*} := \tilde{\mu}$, that is, $\tilde{\mu}^{0,*}$ is the "prior" CRM and $\nu^{0,*}$ is its intensity measure. For every $n \geq 0$, $q, p \geq 1$, we denote by

$$L^p((\nu^{n,*})^q) = L^p((\mathbb{R}^+ \times \mathbb{X})^q, (\mathscr{B}(\mathbb{R}^+) \otimes \mathscr{X})^q, (\nu^{n,*})^q)$$

the Banach space of real-valued functions $f$ on $(\mathbb{R}^+ \times \mathbb{X})^q$, such that $|f|^p$ is integrable with respect to $(\nu^{n,*})^q := (\nu^{n,*})^{\otimes q}$. We write $L^p((\nu^{n,*})^1) = L^p(\nu^{n,*})$, $p \geq 1$. The symbol $L_s^2((\nu^{n,*})^2)$ is used to denote the Hilbert subspace of $L^2((\nu^{n,*})^2)$ generated by the *symmetric functions* on $(\mathbb{R}^+ \times \mathbb{X})^2$. Note that a function $f$, on $(\mathbb{R}^+ \times \mathbb{X})^2$, is said to be symmetric whenever $f(s, x; t, y) = f(t, y; s, x)$ for every $(s, x), (t, y) \in \mathbb{R}^+ \times \mathbb{X}$.

Now we introduce various kernels which will enter either the statements or the conditions of the posterior CLTs. For $n \geq 0$, we denote the posterior hazard rate and posterior cumulative hazard, given $\mathbf{X}$ and $\mathbf{Y}$, by

$$\tag{18} \tilde{h}_{\Delta^{n,*}}(t) = \int_{\mathbb{X}} k(t, x)[\tilde{\mu}^{n,*}(dx) + \Delta^{n,*}(dx)] = \tilde{h}^{n,*}(t) + \sum_{i=1}^{k} J_i k(t, X_i^*)$$

$$\tag{19} \tilde{H}_{\Delta^{n,*}}(T) = \int_0^T \tilde{h}_{\Delta^{n,*}}(t) \, dt = \tilde{H}^{n,*}(T) + \sum_{i=1}^{k} J_i \int_0^T k(t, X_i^*) \, dt.$$

In (18) and (19), we implicitly introduced the notation $\tilde{h}^{n,*}(t)$ and $\tilde{H}^{n,*}(T)$ for, respectively, the hazard rate and cumulative hazard without fixed points of discontinuity. Note that $\tilde{h}_{\Delta^{0,*}}(t)$ coincides with $\tilde{h}(t)$, the prior hazard rate.

Furthermore, we need to define two basic classes of kernels:

(i) for every $n \geq 0$ and every $f \in L_s^2((\nu^{n,*})^2)$, the kernel $f \star_{1,n}^1 f$ is defined on $(\mathbb{R}^+ \times \mathbb{X})^2$ and is equal to the *contraction*

$$(20) \qquad f \star_{1,n}^1 f(t_1, x_1; t_2, x_2) = \int_{\mathbb{R}^+ \times \mathbb{X}} f(t_1, x_1; s, y) f(s, y; t_2, x_2) \nu^{n,*}(ds, dy);$$

(ii) for every $n \geq 0$ and every $f \in L_s^2((\nu^{n,*})^2)$, the kernel $f \star_{2,n}^1 f$ is defined on $(\mathbb{R}^+ \times \mathbb{X})$ and is given by

$$(21) \qquad f \star_{2,n}^1 f(t, x) = \int_{\mathbb{R}^+ \times \mathbb{X}} f(t, x; s, y)^2 \nu^{n,*}(ds, dy).$$

The "star" notation is rather common, see, for example, [16, 28, 34]. Note that the Cauchy–Schwarz inequality yields that $f \star_{1,n}^1 f \in L_s^2((\nu^{n,*})^2)$. It is worth noting that the two operators "$\star_{1,n}^1$" and "$\star_{2,n}^1$," which appear in the statements of our CLTs, can be used to obtain explicit (combinatorial) expressions of the moments and of the cumulants associated with single and double integrals with respect to a Poisson (completely) random measure. See [31] for a discussion of this point.

Introduce now a last set of kernels which will appear in the conditions of the results discussed in Section 4. Fix $n \geq 0$, take $T$ such that $0 \leq T < +\infty$ and define

$$(22) \qquad k_T^{(0)}(s, x) = s \int_0^T k(t, x) \, dt, \qquad (s, x) \in \mathbb{R}^+ \times \mathbb{X};$$

$$(23) \qquad k_T^{(1)}(s, x; t, y) = \frac{st}{T} \int_0^T k(u, x) k(u, y) \, du;$$

$$(24) \qquad k_T^{(2)}(s, x) = \frac{s^2}{T} \int_0^T k(u, x)^2 \, du;$$

$$(25) \qquad k_{T,n}^{(3)}(s, x) = \int_{\mathbb{R}^+ \times \mathbb{X}} k_T^{(1)}(s, x; u, w) \nu^{n,*}(du, dw).$$

Finally, for $(s, x) \in \mathbb{R}^+ \times \mathbb{X}$ define the *random* kernel

$$
\begin{aligned}
(26) \qquad k_{T, \Delta^{n,*}}^{(4)}(s, x) &= \frac{s}{T} \int_0^T k(u, x) \int_{\mathbb{X}} k(u, y) \Delta^{n,*}(dy) \, du \\
&= \sum_{i=1}^k k_T^{(1)}(s, x; J_i, X_i^*).
\end{aligned}
$$

4.2. *General results.* Before stating the results concerning the asymptotic behavior of functionals of random hazards, we need to make some more technical assumptions, which do not appear to be very restrictive; indeed, in the following examples, involving kernels and CRMs commonly exploited in practice, they will be shown to hold.

In the sequel we consider mixture hazards (1) which, in addition to (H1)–(H2), satisfy also

(H3)
$$\int_{\mathbb{R}^+ \times \mathbb{X}} k(t,x)^j v^j \rho(dv|x)\lambda(dx) < +\infty \qquad \forall t, \ j = 1, 2, 4;$$

$$\int_0^T \int_{\mathbb{R}^+ \times \mathbb{X}} k(t,x)^j v^j \rho(dv|x)\lambda(dx)\,dt < +\infty \qquad \forall T \geq 0, j = 2, 4.$$

See [27, 28] for a discussion of these conditions. Recall from (18), that $\tilde{h}^{n,*}(t)$ stands for the posterior hazard without fixed points of discontinuity (given **X** and **Y**) and is characterized by (13). It is straightforward to see that, if the prior hazard rate satisfies (H1)–(H3), then $\tilde{h}^{n,*}(t)$ meets (H1)–(H3) as well.

Given an event $B \in \mathscr{F}$, we will say that $B$ has $\mathbb{P}\{\cdot|\mathbf{X}, \mathbf{Y}\}$-*probability* 1 whenever there exists $\Omega' \in \mathscr{F}$ such that $\mathbb{P}\{\Omega'\} = 1$, and, for every fixed $\omega \in \Omega'$, the random probabilty measure $A \mapsto \mathbb{P}\{\mathbf{X} \in A|\mathbf{Y}\}(\omega)$ has support contained in the set of those $(x_1, \ldots, x_n) \in \mathbb{X}^n$ such that

$$\mathbb{P}\{B|\mathbf{X} = (x_1, \ldots, x_n), \mathbf{Y}\} = 1.$$

Finally, fix a sample size $n \geq 1$ for the remainder of the section. The following Theorems 8, 9 and 10 provide sufficient conditions to have that linear and quadratic functionals associated with posterior random hazard rates verify a CLT. The first result deals with linear functionals.

THEOREM 8 (Linear functionals). *Suppose:* (i) $k_T^{(0)} \in L^3(\nu^{n,*})$ *for every* $T > 0$; (ii) *there exists a strictly positive (deterministic) function* $T \mapsto C_0(n, k, T)$ *such that, as* $T \to +\infty$,

$$(27) \qquad C_0^2(n, k, T) \times \int_{\mathbb{R}^+ \times \mathbb{X}} [k_T^{(0)}(s, x)]^2 \nu^{n,*}(ds, dx) \to \sigma_0^2(n, k),$$

$$(28) \qquad C_0^3(n, k, T) \times \int_{\mathbb{R}^+ \times \mathbb{X}} [k_T^{(0)}(s, x)]^3 \nu^{n,*}(ds, dx) \to 0,$$

*where* $\sigma_0^2(n, k) \in (0, +\infty)$. *Also assume that, with* $\mathbb{P}\{\cdot|\mathbf{X}, \mathbf{Y}\}$-*probability 1,*

$$(29) \qquad \lim_{T \to +\infty} C_0(n, k, T) \times \sum_{i=1}^k J_i \int_0^T k(t, X_i^*)\,dt = m(n, \Delta^{n,*}, k) \in [0, +\infty).$$

*Then, a.s.-$\mathbb{P}$, for every real $\lambda$,*

$$\mathbb{E}[\exp(i\lambda C_0(n, k, T)[\tilde{H}(T) - \mathbb{E}[\tilde{H}^{n,*}(T)]])|\mathbf{Y}]$$

$$\xrightarrow[T \to +\infty]{} \mathbb{E}\left[\exp\left(i\lambda m(n, \Delta^{n,*}, k) - \frac{\lambda^2}{2}\sigma_0^2(n, k)\right)\Big|\mathbf{Y}\right].$$

REMARK 2. When $n = 0$ and setting, by convention, $\mathbf{Y} = \mathbf{X} = 0$ so that $\sigma\{\mathbf{Y}, \mathbf{X}\} = \{\Omega, \varnothing\}$, one recovers Theorem 1 in [27] for prior random hazards. The same applies for the following two results concerning path-second moments and path-variances.

THEOREM 9 (Path-second moments). *Suppose $k_{T,n}^{(3)} \in L^2(\nu^{n,*}) \cap L^1(\nu^{n,*})$, $k_T^{(2)} \in L^3(\nu^{n,*})$ and that there exists a strictly positive function $C_1(n, k, T)$ such that the following asymptotic conditions are satisfied as $T \to +\infty$:*

1. $2C_1^2(n, k, T) \|k_T^{(1)}\|_{L^2((\nu^{n,*})^2)}^2 \to \sigma_1^2(n, k) \in (0, +\infty)$;
2. $C_1^4(n, k, T) \|k_T^{(1)}\|_{L^4((\nu^{n,*})^2)}^4 \to 0$;
3. $C_1^4(n, k, T) \|k_T^{(1)} \star_{1,n}^1 k_T^{(1)}\|_{L^2((\nu^{n,*})^2)}^2 \to 0$;
4. $C_1^4(n, k, T) \|k_T^{(1)} \star_{2,n}^1 k_T^{(1)}\|_{L^2(\nu^{n,*})}^2 \to 0$;
5. $C_1^2(n, k, T) \|k_T^{(2)} + 2k_{T,n}^{(3)} + 2k_{T,\Delta^{n,*}}^{(4)}\|_{L^2(\nu^{n,*})}^2 \to \sigma_4^2(n, \Delta^{n,*}, k) \in [0, +\infty)$, *with $\mathbb{P}\{\cdot | \mathbf{X}, \mathbf{Y}\}$-probability 1*;
6. $C_1^3(n, k, T) \|k_T^{(2)} + 2k_{T,n}^{(3)} + 2k_{T,\Delta^{n,*}}^{(4)}\|_{L^3(\nu^{n,*})}^3 \to 0$, *with $\mathbb{P}\{\cdot | \mathbf{X}, \mathbf{Y}\}$-probability 1*;
7. *with $\mathbb{P}\{\cdot | \mathbf{X}, \mathbf{Y}\}$-probability 1,*

$$\frac{C_1(n, k, T)}{T} \int_0^T \left(\sum_{j=i}^k J_i k(t, X_i^*)\right)^2 dt \to v(n, \Delta^{n,*}, k) \in [0, +\infty).$$

*Moreover, define*

$$(30) \qquad A_T^{n,*} := \sum_{j=1}^k \frac{2J_j}{T} \int_0^T \mathbb{E}[\tilde{h}^{n,*}(t)] k(t, X_j^*) \, dt.$$

*Then, a.s.-$\mathbb{P}$, for every real $\lambda$,*

$$\mathbb{E}\left[\exp\left(i\lambda C_1(n, k, T) \left\{\frac{1}{T} \int_0^T \tilde{h}(t)^2 \, dt - A_T^{n,*} - \frac{1}{T} \int_0^T \mathbb{E}[\tilde{h}^{n,*}(t)^2] \, dt\right\}\right) \Big| \mathbf{Y}\right]$$

$$\xrightarrow[T \to +\infty]{} \mathbb{E}\left[\exp\left(i\lambda v(\Delta^{n,*}, k) - \frac{\lambda^2}{2}(\sigma_1^2(n, k) + \sigma_4^2(n, \Delta^{n,*}, k))\right) \Big| \mathbf{Y}\right].$$

THEOREM 10 (Path-variances). *Suppose that the assumptions of Theorem 8 and Theorem 9 are satisfied. Assume, moreover, that*

1. $C_1(n, k, T)/(TC_0(n, k, T))^2 \to 0$;
2. $2C_1(n, k, T)\mathbb{E}[\tilde{H}^{n,*}(T)]/(T^2 C_0(n, k, T)) \to \delta(n, k) \in \mathbb{R}$;
3. $\|C_1(n, k, T)(k_T^{(2)} + 2k_{T,n}^{(3)} + 2k_{T,\Delta^{n,*}}^{(4)}) - \delta(n, k)C_0(n, k, T)k_T^{(0)}\|_{L^2(\nu^{n,*})}^2 \to \sigma_5^2(n, \Delta^{n,*}, k) \in [0, +\infty)$, *with $\mathbb{P}\{\cdot | \mathbf{X}, \mathbf{Y}\}$-probability 1*

*and $A_T^{n,*}$ is given by (30). Then, a.s.-$\mathbb{P}$, for every real $\lambda$*

$$\mathbb{E}[e^{i\lambda C_1(n,k,T)\{1/T\int_0^T[\tilde{h}(t)-\tilde{H}(T)/T]^2\,dt-A_T^{n,*}-1/T\int_0^T \mathbb{E}[\tilde{h}^{n,*}(t)^2]\,dt+\mathbb{E}[\tilde{H}^{n,*}(T)]^2/T^2\}}|\mathbf{Y}]$$

$$\underset{T\to+\infty}{\longrightarrow} \mathbb{E}[e^{i\lambda(v(n,\Delta^{n,*},k)-\delta(n,k)m(n,\Delta^{n,*},k))-\lambda^2/2(\sigma_1^2(n,k)+\sigma_5^2(n,\Delta^{n,*},k))}|\mathbf{Y}].$$

REMARK 3. We stress that, in general, the four quantities $m(n,\Delta^{n,*},k)$, $\sigma_4^2(n,\Delta^{n,*},k)$, $v(n,\Delta^{n,*},k)$ and $\sigma_5^2(n,\Delta^{n,*},k)$ (appearing in the previous three statements) can be random.

4.3. *Applications.* In this section we derive CLTs for functionals of posterior hazards based on the four kernels (2)–(5), combined with generalized gamma CRMs [2], namely CRMs as in (7) with $\gamma$ a positive constant. The measure $\lambda$ is chosen such that the life-testing model is well defined and (H1)–(H3) are met. Many other classes of CRM represent possible alternatives and one can proceed as below. It is important to recall that consistency of all the models dealt with below is easily deduced from the results in Section 3.

In all the cases we get to the conclusion that the asymptotic behavior of functionals of the posterior hazard rate coincides exactly with the behavior of functionals of the prior hazard. To see why this happens, let us focus on the behavior of the trend of the posterior CRM. It turns out that $\mathbb{E}[\tilde{H}^{n,*}(T)] \sim \psi_1(T) + \psi_2(T;\mathbf{Y})$, where $\psi_1(T) \sim \mathbb{E}[\tilde{H}(T)]$ and $\psi_2(T;\mathbf{Y})$ explicitly depends on the data $\mathbf{Y}$, is different from 0 for every $T > 0$ and $\psi_2(T;\mathbf{Y}) = o(\psi_1(T))$. Moreover, once the rate of divergence from the trend $C_0(n,k,T)$ is computed, one finds that $C_0(n,k,T) = C_0(k;T)$ and $C_0(k,T)^{-1} \times \psi_2(T;\mathbf{Y}) \to 0$ as $T \to \infty$. To fix ideas, consider a DL mixture hazard with generalized gamma CRM given one observation $Y_1$: one obtains

$$\mathbb{E}[\tilde{H}^{1,*}(T)] = \frac{T^2}{2\gamma^{1-\sigma}} - T\left[\frac{Y_1}{\gamma^{1-\sigma}} - \int_0^{Y_1}\frac{1}{(Y_1-x+\gamma)^{1-\sigma}}\,dx\right] + O(1)$$

and, since the divergence rate $C_0(n,k,T)^{-1}$ is equal to $T^{3/2}$, the influence of the data vanishes at a rate $T^{-1/2}$. Similar phenomena occur when studying the asymptotic behavior of the part of the posterior corresponding to the fixed points of discontinuity. This basically explains why the forthcoming CLTs do not depend on the data. Such an outcome is quite surprising, at least to us. Note, indeed, that the Poisson intensity of the posterior CRM (13) depends explicitly on the data $\mathbf{Y}$, which implies that the posterior hazard, and a fortiori the posterior cumulative hazard, depend on the data for any $T$. Also, the fact that the variance of the asymptotic Gaussian distribution is not influenced by the data is somehow counterintuitive: since the contribution of the CRM vanishes in the limit, one would expect the variance to become smaller and smaller as more data come in. Since this does not happen, our findings provide some evidence that the choice of the CRM really matters whatever the size of the dataset. Hence, one should carefully select the kernel and CRM so to incorporate prior knowledge appropriately into the model; the neat CLTs presented here provide a guideline in this respect by highlighting trend, oscillation around the trend and asymptotic variance.

### 4.3.1. *Asymptotics for kernels with finite support.*

We start by considering kernels with finite support, namely, the DL, OU and rectangular ones with generalized gamma CRM and take $\lambda$ to be the Lebesgue measure on $\mathbb{R}^+$. This ensures that (H1)–(H3) are satisfied. For a generalized gamma CRM one has, for any $c > 0$, $\int_0^\infty s^c \rho(ds) = [(1-\sigma)_{c-1}] \gamma^{-c+\sigma} := K_\rho^{(c)}$, where $(a)_n := \Gamma(a+n)/\Gamma(a)$ denotes the Pochhammer symbol. Since in the posterior the CRM becomes nonhomogeneous with updated intensity (13), the verification of the conditions of Theorems 8–10 can become cumbersome. However, for any $A \in \mathbb{R}_+^2$, one has

$$\overline{\nu}(A) \leq \nu^{n,*}(A) \leq \nu(A), \tag{31}$$

where $\overline{\nu}(dv, dx) := \exp\{-n k_{Y_{(n)}}^{(0)}(v,x)\} \nu(dv, dx)$ and $Y_{(n)}$ stands for the largest lifetime. Having a lower and an upper bound for the Poisson intensity $\nu^{n,*}$ allows then to use, conditionally on $\mathbf{X}, \mathbf{Y}$, a comparison result analogous to Theorem 4 of [27] in order to check the conditions of the posterior CLTs.

Let us first consider linear functionals for the OU kernel. Note that $k_T^{(0)}(v, x) = v\sqrt{2/\kappa}(1 - e^{-\kappa(T-x)})\mathbb{I}_{(0 \leq x \leq T)}$, and that $k_T^{(0)} \in L^3(\nu)$, so that condition (i) of Theorem 8 is a direct consequence of (31). Next, one can check that $\|k_T^{(0)}\|_{L^2(\overline{\nu})}^2 \sim \|k_T^{(0)}\|_{L^2(\nu)}^2 \sim 2\kappa^{-1} K_\rho^{(2)} T$. In fact, the dominating term in the norm with respect to $\overline{\nu}$ is the integral over $\mathbb{R}^+ \times [Y_{(n)}, \infty)$, which is in turn equal to the dominating term of $\|k_T^{(0)}\|_{L^2(\nu)}^2$. Moreover, $T^{-3/2}\|k_T^{(0)}\|_{L^3(\nu^{n,*})}^3 \to 0$ and we have that (27) and (28) are satisfied with $C_0(n, k, T) = C_0(0, k, T) = 1/\sqrt{T}$ and $\sigma_0^2(n, k) = \sigma_0^2(0, k) = 2\kappa^{-1} K_\rho^{(2)}$, which importantly does not depend on the observations $\mathbf{Y}$. As for (29), with $\mathbb{P}\{\cdot | \mathbf{X}, \mathbf{Y}\}$-probability 1, we have $\sum_{i=1}^k k_T^{(0)}(J_i, X_i^*) = O(T^{-1})$ as $T \to \infty$, so that (29) holds with $m(n, \Delta^{n,*}, k) = 0$, not depending on $\mathbf{X}, \mathbf{Y}$. Finally, since $\|k_T^{(0)}\|_{L^1(\overline{\nu})} \sim \|k_T^{(0)}\|_{L^1(\nu)} \sim K_\rho^{(1)}\sqrt{2/\kappa}T$, then $\mathbb{E}[\tilde{H}^{n,*}(T)] \sim \mathbb{E}[\tilde{H}(T)]$. Hence, from Theorem 8 combined with fact that the limiting mean does not depend on $\mathbf{Y}$, it follows that

$$\mathbb{E}\left[\exp\left(i\lambda \frac{[\tilde{H}(T) - \sqrt{2/\kappa}\gamma^{-1+\sigma}T]}{T^{1/2}}\right) \Big| \mathbf{Y}\right] \xrightarrow[T \to +\infty]{} \exp\left(-\frac{\lambda^2}{2}\sigma_0^2(0, k)\right), \tag{32}$$

where $\sigma_0^2(0, k) = 2\kappa^{-1}(1-\sigma)\gamma^{-2+\sigma}$. Therefore, the posterior cumulative hazard has the same asymptotic behavior as the prior cumulative hazard. As mentioned before, this is quite surprising also in the light of the consistency result.

Let us now consider the path-second moment. We obtain $k_T^{(1)}(v, x; u, y) = \frac{uv}{T} e^{\kappa(x+y)} (e^{-2\kappa(x \vee y)} - e^{-2\kappa T}) \mathbb{I}_{(0 \le x, y \le T)}$, and, as for condition 1 in Theorem 9, one finds that $\|k_T^{(1)}\|_{L^2(\overline{\nu}^2)}^2 \sim \|k_T^{(1)}\|_{L^2(\nu^2)}^2 \sim 2\kappa^{-1}(K_\rho^{(2)})^2 T$ and, hence, $C_1(n, k, T) = \sqrt{T}$ and $\sigma_1^2(n, k) = 2\kappa^{-1}(K_\rho^{(2)})^2$, which coincide with the case $n = 0$. The idea here is the same as before, namely that the dominating term of the norm with respect to $\overline{\nu}^2$ is the integral over $(\mathbb{R}^+ \times [Y_{(n)}, \infty])^2$, which is equal to the dominating term of $\|k_T^{(1)}\|_{L^2(\nu^2)}^2$. Then, conditions 2., 3. and 4. are verified since they are verified for $n = 0$. In particular, note that $k_T^{(1)} \star_{i,n}^1 k_T^{(1)} \le k_T^{(1)} \star_{i,0}^1 k_T^{(1)}$ for $i = 1, 2$. As for condition 5., one first check that $k_{T, \Delta^{n,*}}^{(4)} = O(T^{-1})$, then some tedious algebra allows to verify that it is satisfied with $\sigma_4^2(n, \Delta^{n,*}, k) = K_\rho^{(4)} + \frac{8}{\kappa} K_\rho^{(3)} K_\rho^{(1)} + \frac{16}{\kappa^2} K_\rho^{(2)} (K_\rho^{(1)})^2$. This is, indeed, a delicate point since both $k_{T,n}^{(3)}$ and the norm with respect to the updated Poisson intensity $\nu^{n,*}$ depend on the posterior. Once this is done, it is not difficult to check that condition 6. is satisfied. Moreover, the quantity $v(n, \Delta^{n,*}, k)$ in condition 7. can be shown to be 0, whereas $A_T^{n,*} = O(T^{-1})$ in (30). Finally, one can check that $\frac{1}{T} \int_0^T \mathbb{E}[\tilde{h}^{n,*}(t)^2] \, dt \sim \frac{1}{T} \int_0^T \mathbb{E}[\tilde{h}^{0,*}(t)^2] \, dt \sim K_\rho^{(2)} + \frac{2}{\kappa} (K_\rho^{(1)})^2$, so that, from Theorem 9, we deduce the following CLT for the path-second moment:

$$(33) \quad \mathbb{E}\left[ \exp\left( i\lambda \sqrt{T} \left\{ \frac{1}{T} \int_0^T \tilde{h}(t)^2 \, dt - \gamma^{-2+\sigma} \left( 1 - \sigma + \frac{2\gamma^\sigma}{\kappa} \right) \right\} \right) \Big| \mathbf{Y} \right]$$
$$\underset{T \to +\infty}{\longrightarrow} \exp\left( -\frac{\lambda^2}{2} (\sigma_1^2(n, k) + \sigma_4^2(n, \Delta^{n,*}, k)) \right),$$

where $\sigma_1^2(n, k) + \sigma_4^2(n, \Delta^{n,*}, k) = \frac{(1-\sigma)(16\kappa^{-1}\gamma^{2\sigma} + 2(9-5\sigma)\gamma^\sigma + \kappa(2-\sigma)_2)}{\kappa\gamma^{4-\sigma}}$.

As far as the path-variance is concerned, one verifies easily that the conditions of Theorem 10 are satisfied, with $\delta(n, k) = \frac{2^{3/2}}{\sqrt{\kappa}} K_\rho^{(1)}$ and $\sigma_5^2(n, \Delta^{n,*}, k) = K_\rho^{(4)}$, which again do not depend on the observations $\mathbf{Y}$. As for the posterior mean of the path-variance, one finds that $\frac{1}{T} \int_0^T \mathbb{E}[\tilde{h}^{n,*}(t) - \frac{\mathbb{E}[\tilde{H}^{n,*}(T)]}{T}]^2 \, dt = K_\rho^{(2)} + o(T^{-1/2})$, so that Theorem 10 leads to

$$(34) \quad \mathbb{E}\left[ \exp\left( i\lambda \sqrt{T} \left\{ \frac{1}{T} \int_0^T \left[ \tilde{h}(t) - \frac{\tilde{H}(T)}{T} \right]^2 \, dt - \frac{1-\sigma}{\gamma^{2-\sigma}} \right\} \right) \Big| \mathbf{Y} \right]$$
$$\underset{T \to +\infty}{\longrightarrow} \exp\left( -\frac{\lambda^2}{2} (\sigma_1^2(n, k) + \sigma_5^2(n, \Delta^{n,*}, k)) \right),$$

where $\sigma_1^2(n, k) + \sigma_5^2(n, \Delta^{n,*}, k) = \frac{(1-\sigma)(2(1-\sigma)\gamma^\sigma + \kappa(2-\sigma)_2)}{\kappa\gamma^{4-\sigma}}$.

For the other two kernels, namely rectangular and DL, one can proceed along the same lines of reasoning and, again, the asymptotic posterior behavior coincides with the one of the prior. In particular, one obtains that for linear functionals and quadratical functionals of hazard rates based on the rectangular kernel, the CLTs (32), (33) and (34) hold with the same rate functions and appropriately modified constants and variances (for the exact values see [27], since they coincide with the a priori ones). As for the DL kernel the CLT for the posterior cumulative hazard is of the form

$$\frac{1}{T^{3/2}}\left[\tilde{H}(T) - \frac{1}{2\gamma^{1-\sigma}}T^2\right]\Big|\mathbf{Y} \xrightarrow{\text{law}} X \sim \mathcal{N}\left(0, \frac{1-\sigma}{3\gamma^{2-\sigma}}\right).$$

With reference to quadratic functionals, in this case, some of the conditions of Theorems 9 and 10 are violated already in prior (see [27] for details).

### 4.3.2. *Asymptotics for exponential kernel.*

Here we consider random hazards based on the exponential kernel. Indeed, it is crucial to consider also a kernel with full support, since one may think that the lack of dependence on the data of posterior functionals may be due to the boundedness of the support of the kernels dealt with in Section 4.3.1. However, it turns out that, again, the posterior CLTs coincide with the corresponding prior CLTs.

In particular, set, within (7), $\lambda(dx) = x^{-1/2}e^{-1/x}(2\sqrt{\pi})^{-1}$: this implies that $\tilde{h}(0) < \infty$ a.s., (8) is in order and (H1)–(H3) are satisfied. This model is of interest also beyond the scope of the present asymptotic analysis; in fact, it leads to a prior mean $\mathbb{E}[\tilde{h}(t)] = K_\rho^{(1)}(t+1)^{-1/2}$ and, thus, we have a nonparametric prior centered on a quasi Weibull hazard, which is a desirable feature in survival analysis.

We start by investigating the linear functional of $\tilde{h}$: here we provide details also for the derivation of the prior CLTs since this model has not been considered in [27]. In this case, we have that $k_T^{(0)}(v,x) = v(1 - e^{-T/x})$ and $k_T^{(0)}(v,x) \in L^3(\nu)$ for all $T > 0$ and the same holds for the posterior. We also have that $\|k_T^{(0)}\|_{L^1(\nu)} = K_\rho^{(1)}(\sqrt{1+T} - 1)$, so that, as $T \to \infty$, $\mathbb{E}[\tilde{H}(T)] \sim K_\rho^{(1)}\sqrt{T}$. When considering the posterior, one can check that $\|k_T^{(0)}\|_{L^1(\overline{\nu})} \sim \|k_T^{(0)}\|_{L^1(\nu)}$. In fact, by a change of variable and dominated convergence

$$\frac{\|k_T^{(0)}\|_{L^1(\overline{\nu})}}{\sqrt{T}} = \frac{1}{\sqrt{T}}\int_{\mathbb{R}^+}\frac{(1 - e^{-T/x})}{[\gamma + n(1 - e^{-Y_{(n)}/x})]^{1-\sigma}}\frac{e^{-1/x}x^{-1/2}}{2\sqrt{\pi}}\,dx$$

$$= \frac{1}{2\sqrt{\pi}}\int_{\mathbb{R}^+}\frac{(1 - e^{-y})e^{-y/T}y^{-3/2}}{[\gamma + n(1 - e^{-yY_{(n)}/T})]^{1-\sigma}}\,dy$$

$$\xrightarrow[T \to +\infty]{}\frac{1}{2\sqrt{\pi}}\int_{\mathbb{R}^+}\frac{1 - e^{-y}}{y^{3/2}\gamma^{1-\sigma}}\,dy = K_\rho^{(1)}.$$

Therefore, $\mathbb{E}[\tilde{H}^{n,*}(T)] \sim \mathbb{E}[\tilde{H}(T)]$. Similar arguments lead to show that $\|k_T^{(0)}\|_{L^2(\overline{\nu})}^2 \sim \|k_T^{(0)}\|_{L^2(\nu)}^2 \sim (2-\sqrt{2})K_\rho^{(2)}\sqrt{T}$ and, hence, we have $C_0(n,k,T) = C_0(0,k,T) = T^{-1/4}$ and $\sigma_0^2(n,k) = \sigma_0^2(0,k) = (2-\sqrt{2})K_\rho^{(2)}$. Moreover, $\|k_T^{(0)}\|_{L^3(\nu)}^3 \sim O(\sqrt{T})$ is sufficient for concluding that (28) holds both for the prior and the posterior. Finally, as $T \to \infty$, $\sum_{i=1}^k k_T^{(0)}(J_i, X_i^*) = O(1)$ with $\mathbb{P}\{\cdot|\mathbf{X}, \mathbf{Y}\}$-probability 1; thus, also in this case (29) holds with $m(n, \Delta^{n,*}, k) = 0$. We can then deduce from Theorem 8 that

$$\mathbb{E}\left[\exp\left(i\lambda\frac{[\tilde{H}(T) - \gamma^{-(1-\sigma)}T^{1/2}]}{T^{1/4}}\right)\Big|\mathbf{Y}\right] \underset{T \to +\infty}{\longrightarrow} \exp\left(-\frac{\lambda^2}{2}\sigma_0^2(0,k)\right)$$

for any sample size $n \geq 0$ and with $\sigma_0^2 = (2-\sqrt{2})(1-\sigma)\gamma^{-1+\sigma}$. Hence, we have shown that the exponential kernel hazard exhibits both trend and oscillations of order $T^{1/2}$ and verifies exactly the same CLT for both prior and posterior cumulative hazard, thus confirming that the asymptotics is not influenced by the data.

Our results for quadratic functionals do not apply to the exponential kernel. To see this, note that $k_T^{(1)}(v,x;u,y) = \frac{uv}{t}\frac{1}{x+y}(1 - \exp\{-\frac{x+y}{xy}T\})$ and, by calculating the norm with respect to $\nu^2$, we get

$$\|k_T^{(1)}\|_{L^2(\nu^2)}^2 = \frac{(K_\rho^{(2)})^2}{16(2T^2 + 3T + 1)},$$

which implies $C_1(0,k,T) = T$. However, $\|k_T^{(1)}\|_{L^4(\nu^2)}^4 \sim \frac{d}{T^4}$, $d$ being a positive constant, so that condition 2 in Theorem 9 does not hold.

**5. Concluding remarks.** In the present paper we have investigated two different asymptotic aspects of a random hazard model, namely consistency and the behavior of a functionals of the hazard as time diverges. As for the former, we have provided a general weak consistency criterion for mixture random hazards and established weak consistency for specific models with respect to large classes of "true hazards" $h_0$. It seems worth discussing briefly the case of Weibull hazards, that is, $h_0(t) = \alpha\lambda t^{\alpha-1}(\alpha/\lambda)(t/\lambda)^{\alpha-1}$ with $\alpha, \lambda > 0$, which are widely used in the parametric setup. The case of $\alpha > 1$ is covered by both Theorem 4 (DL kernel) and Theorem 6 (OU kernel). When $\alpha < 1$, $h_0$ is a completely monotone function and it would naturally belong to the domain of attraction of Theorem 7; however, in such a case $h_0(0)$ is not finite and, hence, the required conditions are not met. Nonetheless, $h_0$ can be approximated to any order of accuracy by $h_\varepsilon(t) = (\alpha/\lambda)((t+\varepsilon)/\lambda)^{\alpha-1}$, for some small enough $\varepsilon > 0$, when accuracy is measured in terms of survival functions. In fact, it is easy to see that for $S_0(t)$ and $S_\varepsilon(t)$, the survival functions corresponding to $h_0$ and $h_\varepsilon$, respectively, $\sup_t |S_0(t) - S_\varepsilon(t)|$ goes to zero as $\varepsilon$ approaches zero. Finally, note that Theorem 7 applies to $h_\varepsilon$ for any $\varepsilon > 0$. Further work is needed in order to extend the consistency result to completely monotone hazards which explode in zero; for such cases, condition (17) is probably too strong.

Future work will also focus on achieving consistency with respect to stronger topologies; two are the possible routes in this direction. The first one is to investigate under which additional conditions on the CRM $\tilde{\mu}$ and restrictions on the form of the true hazard rate $h_0$ we get $L_1$-consistency at the density level, that is, (10) with $A_\varepsilon$ being a $L_1$ neighborhood of $f_0$. To this end, one has then to consider the metric entropy of the subset of $\mathbb{F}$ corresponding to the qualitative condition given on $h_0$. Moreover, one has to investigate in detail the support of the prior $\Pi$ on $\mathbb{F}$ via the mapping $\tilde{h} \to \tilde{f} = \tilde{h} \exp(-\int_0^t \tilde{h})$. This appears to be a rather difficult problem because of $\tilde{h}$ appearing twice, and existing results on random mixing densities are not easily extensible. The second strategy consists of investigating consistency directly at the hazard level. Indeed, weak consistency at the density level implies pointwise consistency of the cumulative hazard:

$$\Pi_n \left\{ h : \left| \int_0^T h(t)\,dt - \int_0^T h_0(t)\,dt \right| \le \varepsilon \right\} \to 1 \qquad \text{a.s.-}P_0^\infty$$

for any $\varepsilon, T > 0$. Among stronger topologies, a promising one seems to be the one induced by $\int_0^\infty |h(t) - h_0(t)| S_0(t)\,dt$, where $S_0(t) = \exp\{-\int_0^t h_0(s)\,ds\}$.

With reference to the study of the asymptotic behavior of functionals of the random hazard, a further interesting development consists in studying the joint limit as both the number of observations and time diverge. To achieve such a result, one probably needs to find a right balance in the simultaneous divergence of the sample size and time, which lets the influence of the data emerge.

## APPENDIX: BACKGROUND, ANCILLARY RESULTS AND PROOFS

### A.1. Completely random measures.
Here we highlight some basic facts on CRMs. The reader is referred to [3] and [22] for exhaustive accounts. Consider a measure space $(\mathbb{X}, \mathscr{X})$, where $\mathbb{X}$ is a complete and separable metric space and $\mathscr{X}$ is the usual Borel $\sigma$-field. Introduce a *Poisson random measure* $\tilde{N}$, defined on some probability space $(\Omega, \mathscr{F}, \mathbb{P})$ and taking values in the set of nonnegative counting measures on $(\mathbb{R}^+ \times \mathbb{X}, \mathscr{B}(\mathbb{R}^+) \otimes \mathscr{X})$, with *intensity measure* $\nu$, that is, $\mathbb{E}[\tilde{N}(dv, dx)] = \nu(dv, dx)$ and, for any $A \in \mathscr{B}(\mathbb{R}^+) \otimes \mathscr{X}$ such that $\nu(A) < \infty$, $\tilde{N}(A)$ is a Poisson random variable of parameter $\nu(A)$. Given any finite collection of pairwise disjoint sets, $A_1, \ldots, A_k$, in $\mathscr{B}(\mathbb{R}^+) \otimes \mathscr{X}$, the random variables $\tilde{N}(A_1), \ldots, \tilde{N}(A_k)$ are mutually independent. Moreover, the intensity measure $\nu$ must satisfy $\int_{\mathbb{R}^+} (v \wedge 1) \nu(dv, \mathbb{X}) < \infty$ where $a \wedge b = \min\{a, b\}$.

Let now $(\mathbb{M}, \mathscr{B}(\mathbb{M}))$ be the space of boundedly finite measures on $(\mathbb{X}, \mathscr{X})$, where $\mu$ is said boundedly finite if $\mu(A) < +\infty$ for every bounded measurable set $A$. We suppose that $\mathbb{M}$ is equipped with the topology of vague convergence and that $\mathscr{B}(\mathbb{M})$ is the corresponding Borel $\sigma$-field. Let $\tilde{\mu}$ be a random element, defined on $(\Omega, \mathscr{F}, \mathbb{P})$ and with values in $(\mathbb{M}, \mathscr{B}(\mathbb{M}))$, and suppose that $\tilde{\mu}$ can be represented as a linear functional of the Poisson random measure $\tilde{N}$ (with intensity $\nu$) as $\tilde{\mu}(B) = \int_{\mathbb{R}^+ \times B} s\tilde{N}(ds, dx)$ for any $B \in \mathscr{X}$. From the properties of $\tilde{N}$ it easily follows that $\tilde{\mu}$ is a CRM on $\mathbb{X}$ [21], that is: (i) $\tilde{\mu}(\varnothing) = 0$ a.s.-$\mathbb{P}$; (ii) for any collection of disjoint sets in $\mathscr{X}$, $B_1, B_2, \ldots$, the random variables $\tilde{\mu}(B_1), \tilde{\mu}(B_2), \ldots$ are mutually independent and $\tilde{\mu}(\bigcup_{j\geq 1} B_j) = \sum_{j\geq 1} \tilde{\mu}(B_j)$ holds true a.s.-$\mathbb{P}$.

Now let $\mathscr{G}_\nu$ be the space of functions $g : \mathbb{X} \to \mathbb{R}^+$ such that $\int_{\mathbb{R}^+ \times \mathbb{X}} [1 - e^{-sg(x)}] \nu(ds, dx) < \infty$. Then, the law of $\tilde{\mu}$ is uniquely characterized by its *Laplace functional* which, for any $g$ in $\mathscr{G}_\nu$, is given by

$$(35) \qquad \mathbb{E}[e^{-\int_{\mathbb{X}} g(x)\tilde{\mu}(dx)}] = \exp\left\{ -\int_{\mathbb{R}^+ \times \mathbb{X}} [1 - e^{-sg(x)}]\nu(ds, dx) \right\}.$$

From (35) it is apparent that the law of the CRM $\tilde{\mu}$ is completely determined by the corresponding intensity measure $\nu$. Letting $\lambda$ be a $\sigma$-finite measure on $\mathbb{X}$, we can always write the Poisson intensity $\nu$ as (6), where $\rho : \mathscr{B}(\mathbb{R}^+) \times \mathbb{X} \to \mathbb{R}^+$ is a kernel [i.e., $x \mapsto \rho(C|x)$ is $\mathscr{X}$-measurable for any $C \in \mathscr{B}(\mathbb{R}^+)$ and $\rho(\cdot|x)$ is a $\sigma$-finite measure on $\mathscr{B}(\mathbb{R}^+)$ for any $x$ in $\mathbb{X}$]. Note that the kernel $\rho(dv|x)$ is uniquely determined outside some set of $\lambda$-measure 0, and that such a disintegration is guaranteed by Theorem 15.3.3 in [17]. Finally, recall (see, e.g., Proposition 1 in [30]) that a linear functional of a CRM, $\int_{\mathbb{X}} f(x)\tilde{\mu}(dx)$, is a.s. finite if and only if

$$(36) \qquad \int_{\mathbb{R}^+ \times \mathbb{X}} [1 - e^{-u|f(x)|v}]\rho(dv|x)\lambda(dx) < +\infty \qquad \forall u > 0.$$

### A.2. Proofs of the results of Section 3.

*Proof of Theorem 2.* The first step consists in adapting the K–L condition (15) to the case of right-censoring. Denote by $\mathbb{F}_0 \subset \mathbb{F} \times \mathbb{F}$ the class of all pairs of density functions $(f_1, f_2)$ such that both $f_1$ and $f_2$ are supported on the entire positive real line. Let $X_i \sim f_i$, for $i = 1, 2$, suppose $X_1$ is stochastically independent of $X_2$ and define $\psi(X_1, X_2) = (X_1 \wedge X_2, \mathbb{I}_{(X_1 \leq X_2)})$. The density of $\psi$ with respect to the Lebesgue measure and the counting measure on $\{0, 1\}$ is given by

$$\phi(f_1, f_2)(z, 1) = f_1(z)\int_z^\infty f_2(x)\,dx, \qquad \phi(f_1, f_2)(z, 0) = \int_z^\infty f_1(x)\,dx f_2(z).$$

Then $\phi$ is one-to-one on $\mathbb{F}_0$ and the maps $\phi$, $\phi^{-1}$ defined on $\mathbb{F}_0$ and $\mathbb{F}_0^* = \phi(\mathbb{F}_0)$, respectively, are continuous with respect to the supremum distance on distribution functions. See Peterson [29]. Denote by $\bar{\Pi}$ the prior on $\mathbb{F}_0$ and by $\Pi^* = \bar{\Pi} \circ \phi^{-1}$ the induced prior on $\mathbb{F}_0^*$. Since $(f_0, f_c) \in \mathbb{F}_0$ by hypothesis, the continuity of $\phi^{-1}$ implies that the posterior $\bar{\Pi}(\cdot|(Z_1, \Delta_1), \ldots, (Z_n, \Delta_n))$ is weakly consistent at $(f_0, f_c)$ if $\Pi^*(\cdot|(Z_1, \Delta_1), \ldots, (Z_n, \Delta_n))$ is weakly consistent at $\phi(f_0, f_c)$. Indicate by $p(x, d)$, for $x \in \mathbb{R}$ and $d = 0, 1$ a generic element of $\mathbb{F}_0^*$. Then, K–L support condition on $\Pi^*$ at $p_0 \in \mathbb{F}_0^*$ takes the form

$$\Pi^*\left\{p : \int_0^\infty p_0(z, 1) \log \frac{p_0(z, 1)}{p(z, 1)} \, dz + \int_0^\infty p_0(z, 0) \log \frac{p_0(z, 0)}{p(z, 0)} \, dz < \varepsilon \right\} > 0$$

for any $\varepsilon > 0$. As observed in Section 2, since the prior on $f_c$ does not play any role in the analysis, we may treat $f_c$ as fixed, that is, take a prior on $\mathbb{F} \times \mathbb{F}$ of the form $\Pi \times \delta_{f_c}$. Hence, by setting $p_0(x, d) = \phi(f_0, f_c)(z, d)$, the K–L condition boils down to

$$\text{(37)} \quad \Pi\left\{f : \int_0^\infty f_0(t) S_c(t) \log \frac{f_0(t)}{f(t)} \, dt + \int_0^\infty S_0(t) f_c(t) \log \frac{S_0(t)}{S_f(t)} \, dt < \varepsilon \right\} > 0$$

for any $\varepsilon > 0$, where we defined the survival functions $S_0(t) = 1 - \int_t^\infty f_0(x) \, dx$, $S_f(t) = 1 - \int_t^\infty f(x) \, dx$ and $S_c(t) = 1 - \int_t^\infty f_c(x) \, dx$.

The next step consists in showing that, under the stated hypotheses, the K–L support condition (37) is satisfied, which in turn implies weak consistency. Specifically, we show that a sufficient condition for (37) is that, for any $\delta > 0$, there exists $T'$ such that, for any $T > T'$,

$$\text{(38)} \quad \Pi\left\{h : \sup_{t \leq T} |h(t) - h_0(t)| < \delta, \int_T^\infty |H - H_0| f_0 < \delta \right\} > 0,$$

where $H(t) = \int_0^t h(s) \, ds$ and $H_0(t) = \int_0^t h_0(s) \, ds$. By the structural properties of the model with (2)–(5), it follows that (38) holds under condition (17) and $\int_0^\infty |\tilde{H}(t) - H_0(t)| f_0(t) \, dt < \infty$ a.s. In particular, the latter is implied by condition (i) and the fact that $\int_0^\infty H_0(t) f_0(t) < \infty$.

Define the set

$$\text{(39)} \quad V(\delta, T) := \left\{h : \sup_{t \leq T} |h(t) - h_0(t)| < \delta, \int_T^\infty |H - H_0| f_0 < \delta \right\},$$

which, by (38), has positive probability for any $\delta$ and any $T$ larger than a time point $T'$ that may depend on $\delta$. Our goal is then to show that, for any $\varepsilon > 0$, there exists $\delta > 0$ and $T$ sufficiently large such that, for any $h \in V(\delta, T)$,

$$\text{(40)} \quad \int_T^\infty \log(f_0/f) f_0 S_c + \log(S_0/S_f) f_c S_0 < \varepsilon/2,$$

$$\text{(41)} \quad \int_0^T \log(f_0/f) f_0 S_c + \log(S_0/S_f) f_c S_0 < \varepsilon/2,$$

where $f(t) = h(t) \exp(-\int_0^t h(s) \, ds)$. Let us start from (40) by noting that

$$
\begin{aligned}
(42) \quad & \int_T^\infty \log\left(\frac{f_0}{f}\right) f_0 S_c + \log\left(\frac{S_0}{S_f}\right) f_c S_0 \\
& \leq \int_T^\infty \log(h_0) f_0 S_c - \int_T^\infty \log(h) f_0 S_c + \int_T^\infty |H - H_0|(f_0 S_c + f_c S_0).
\end{aligned}
$$

As for the first integral in the right-hand side of (42), it is easy to see that $\int_T^\infty \log(h_0) f_0 S_c$ goes to zero as $T \to \infty$. As for the second integral, one needs to consider the case of $h(t)$ that eventually goes to zero, but then the negligibility of the integral as $T \to \infty$ is guaranteed by condition (i) and (8), which is needed for the model to be well defined. As for the third integral in the right-hand side of (42), notice that $f_0(t) S_c(t) + f_c(t) S_0(t) \leq 2 f_0(t)$ for $t$ sufficiently large since $S_c \leq 1$ and $f_c$ is eventually smaller than $h_0$. Therefore $\int_T^\infty |H - H_0|(f_0 S_c + f_c S_0) < 2\delta$ and we can conclude that there exists a positive $\delta$ sufficiently smaller than $\varepsilon/4$ and $T$ sufficiently large such that (40) holds for any $h \in V(\delta, T)$.

We are now left to show that (41) holds. Assume first that $h_0(0) > 0$ and write

$$
\begin{aligned}
(43) \quad & \int_0^T \log(f_0/f) f_0 S_c + \log(S_0/S_f) f_c S_0 \\
& = \int_0^T \log\left(\frac{h_0(t)}{h(t)}\right) f_0(t) S_c(t) \, dt \\
& \quad + \int_0^T \int_0^t [h(s) - h_0(s)] \, ds \, [f_0(t) S_c(t) + f_c(t) S_0(t)] \, dt := I_1 + I_2.
\end{aligned}
$$

Next, let $c := \inf_{t \leq T} h_0(t)$, which is positive by condition (i), and note that, for $\delta < c$ and $h \in V(\delta, T)$,

$$
\begin{aligned}
I_1 & \leq \int_0^T \left|\frac{h_0(t)}{h(t)} - 1\right| f_0(t) S_c(t) \, dt \leq \int_0^T \frac{\delta}{c - \delta} f_0(t) \, dt \leq \frac{\delta}{c - \delta} \\
I_2 & \leq \int_0^T \left[\sup_{s \leq t} |h(s) - h_0(s)|\right] t [f_0(t) S_c(t) + f_c(t) S_0(t)] \, dt \\
& \leq \delta \int_0^\infty t [f_0(t) S_c(t) + f_c(t) S_0(t)] \, dt \leq \delta E_0,
\end{aligned}
$$

where $E_0 := \int_0^\infty t f_0(t) \, dt$ is finite by condition (i) and the last inequality follows from $f_0 S_c + f_c S_0$ being the density of $Z = Y \wedge C$ which, in turn, is stochastically smaller than $Y$. Hence, $I_1 + I_2 \leq \delta(c - \delta)^{-1} + \delta E_0$, so that $\delta < \min\{c\varepsilon/(4+\varepsilon), \varepsilon/(4E_0)\}$ implies (41) for any $h \in V(\delta, T)$, no matter how large $T$ is. Finally, one can choose $\delta$ small enough and $T$ large enough such that (40) and (41) are simultaneously satisfied for any $h \in V(\delta, T)$.

By allowing $h_0(0) = 0$, we need a different bound for $I_1$ in (43). We proceed by taking $0 < \varsigma < T$ and split $I_1$ into

$$I_1 = \int_0^\varsigma \log\left(\frac{h_0(t)}{h(t)}\right) f_0(t) S_c(t) \, dt + \int_\varsigma^T \log\left(\frac{h_0(t)}{h(t)}\right) f_0(t) S_c(t) \, dt := I_{11} + I_{12}.$$

As for $I_{12}$, for fixed $\varepsilon$, find $\delta$ and $T$ such that $h \in V(\delta, T)$ implies $I_{12} + I_2 < \varepsilon/4$, for any $\varsigma$. As for $I_{11}$, we need to prove that, for the same $\varepsilon$ fixed above, there exists a small enough $\varsigma > 0$ such that

$$(44) \qquad \int_0^\varsigma \log(h_0(t)/h(t)) f_0(t) S_c(t) \, dt < \varepsilon/4.$$

This is tantamount of showing that $\log(h_0/\tilde{h}) f_0 S_c$ is integrable in $0$ a.s., which in turn reduces to show that $\log(h_0/\tilde{h}) f_0$ is integrable in $0$ a.s. since $S_c(0) = 1$. Note that it is sufficient to control the worst case, namely when $\tilde{h}(0) = 0$ a.s., but then integrability in $0$ follows from condition (ii). Indeed, we need to show that there exists $0 < p < 1$ such that

$$\limsup_{\tau \downarrow 0} \frac{\log\{h_0(\tau)/\tilde{h}(\tau)\} f_0(\tau)}{\tau^{p-1}} = 0 \qquad \text{a.s.}$$

First note that $\lim_{\tau \downarrow 0} \log\{h_0(\tau)\} f_0(\tau) = 0$. This can be deduced by reasoning in terms of $\log(f_0) f_0$ since, clearly, $h_0(\tau) \sim f_0(\tau)$ as $\tau \to 0$. As for $\log(f_0) f_0$ vanishing at zero, we start considering $f_0$ having regular variation of exponent $0 < p < 1$ at zero, that is, $f_0(\tau) \sim \tau^p L(1/\tau)$ as $\tau \to 0$, for $L(\cdot)$ a slowly varying function at $\infty$. Recall that a positive function $L(x)$ defined on $\mathbb{R}^+$ varies slowly at $\infty$ if, for every fixed $x$, $L(tx)/L(x) \to 1$ as $t \to \infty$. Hence,

$$f_0(\tau) \log[f_0(\tau)] \sim \tau^p \{\log(\tau^p) + \log[L(1/\tau)]\} := \tau^p L^*(1/\tau),$$

where $L^*$ is a slowly varying function at $\infty$. Hence $f_0 \log(f_0)$ is a regularly varying function at zero with exponent $p$ and, in turn, it vanishes in zero. Note that the larger $p$ is, the faster $\log\{h_0(\tau)\} f_0(\tau)$ vanishes as $\tau \to 0$. Next, we have that, for any $0 < p < 1$,

$$\limsup_{\tau \downarrow 0} \frac{\log\{h_0(\tau)/\tilde{h}(\tau)\} f_0(\tau)}{\tau^{p-1}} = 0 + \limsup_{\tau \downarrow 0} \frac{-\log\{\tilde{h}(\tau)\}}{\tau^{p-1}} \leq \lim_{\tau \downarrow 0} \frac{-\log\{\tau^r\}}{\tau^{p-1}},$$

where the last limit is zero for any $0 < p < 1$. The integrability then follows for any $0 < p < 1$. Slightly different arguments can be used when $f_0$ has regular variation of exponent $p > 1$ at zero, while the special case of $f_0$ slowly varying at zero (i.e., $p = 0$) can be dealt by using Lemma 2 of Feller [7], Section VII.8. The proof is then complete.

*Proof of Proposition 3.* The fact that (ii1) and (ii2) are sufficient for condition (ii)(b) of Theorem 2 to hold is straightforward.

Since for DL mixture hazards $\tilde{h}(t) = \tilde{\mu}([0, t])$ and for OU mixtures $\tilde{h}(t) \geq \sqrt{2\kappa}e^{-\kappa\varepsilon}\tilde{\mu}([0, t])$ for any $\varepsilon > t$, condition (ii1) is met for both.

Let us now show that CRMs as in (7) with $\sigma \in (0, 1)$ and $\lambda(dx) = dx$ satisfy condition (ii2). Assume, for the moment, that $\gamma$ in (7) is constant and denote it by $\bar{\gamma}$. Hence, we have the generalized gamma subordinator, whose Laplace exponent is given by $\psi(u) := \sigma^{-1}(u + \bar{\gamma})^{\sigma} - \bar{\gamma}^{\sigma}$. Moreover, $\int_0^{\infty} v^{\varepsilon}\rho(dv) = \infty$ for any $\varepsilon < \sigma$ and the inverse of $\psi(u)$ is of the form $\psi^{-1}(y) = (\sigma y + \gamma^{\sigma})^{1/\sigma} - \gamma$. Thus, we are in a position to apply Proposition 47.18 in [32], which, in our case allow to state that there exists a constant $C$ such that

$$(45) \qquad \liminf_{t \downarrow 0} \frac{\tilde{\mu}([0, t])}{g(t)} = C \qquad \text{a.s. with } 0 < C < \infty,$$

where $g(t) = \log\log(1/t)[(\sigma t^{-1}\log\log(1/t) + \gamma^{\sigma})^{1/\sigma} - \gamma]^{-1}$. From (45) it follows immediately that, for any $\delta > 0$, $\liminf_{t \downarrow 0} \frac{\tilde{\mu}([0, t])}{t^{1/\sigma + \delta}} = \infty$ a.s. Hence, condition (ii2) is satisfied by taking $r = 1/\sigma + \delta$. To see that condition (ii2) holds also if $\tilde{\mu}$ is a nonhomogeneous CRM it is enough to note that the corresponding Laplace exponent $\sigma^{-1}\int_0^{\infty}[(u + \gamma(x))^{\sigma} - \gamma(x)^{\sigma}]\,dx$ is bounded above by $\psi(u) := \sigma^{-1}(u + \bar{\gamma})^{\sigma} - \bar{\gamma}^{\sigma}$ with $\bar{\gamma} = \inf_{x \in \mathbb{R}^+} \gamma(x) \geq 0$ and that, infinitesimally, a nonhomogeneous CRM behaves like a homogeneous one.

*An auxiliary lemma.* Before getting into the proofs of the consistency results, we provide a useful auxiliary result. Let $\mathbb{M}$ be the space of boundedly finite measures on $\mathbb{R}^+$ and denote by $\mathbb{G}$ the space of distribution function associated to it: clearly, any $G \in \mathbb{G}$ will be a nondecreasing càdlàg function on $\mathbb{R}^+$ such that $G(0) = 0$.

LEMMA 11. *Let $\tilde{\mu}$ be a CRM on $\mathbb{R}^+$, satisfying (H1), and denote by $Q$ the distribution induced on $\mathbb{G}$. Then, for any $G_0 \in \mathbb{G}$, any finite $M$ and $\eta > 0$,*

$$Q\left\{G \in \mathbb{G} : \sup_{x \leq M}|G(x) - G_0(x)| < \eta\right\} > 0.$$

PROOF. Fix $\varepsilon > 0$ and choose $(z_0, \ldots, z_N)$ such that (i) $0 = z_0 \leq z_1 < \cdots < z_N = M$; (ii) all locations, where $G_0$ has a jump of size larger than $\varepsilon/2$, are contained in $(z_1, \ldots, z_{N-1})$; (iii) for $l = 1, \ldots, N$, $G_0(z_l^-) - G_0(z_{l-1}) \leq \varepsilon$. Next, define

$$(46) \qquad G_{\varepsilon}(x) = \sum_{l=1}^{N} j_l \mathbb{I}_{(z_l \leq x)},$$

where the jump $j_l$ at $z_l$ is given by $j_l = G_0\{z_l\} + G_0(z_l^-) - G_0(z_{l-1})$, for $l = 1, \ldots, N$. If $z_1 = 0$, then set by convention $G_0(z_0) := G_0(0^-) = 0$ and $G_\varepsilon(z_0) := G_\varepsilon(0^-) = 0$. By construction $G_\varepsilon(x) \leq G_0(x)$ for any $x \leq M$ and $\sup_{x \leq M} [G_0(x) - G_\varepsilon(x)] \leq \varepsilon$. Under (H1), it can be proved that, for any $x \leq M$ and for any $a, b$ such that $0 < a < b$, $Q\{G \in \mathbb{G} : G(x) \in (a, b)\} > 0$. See, for example, Proposition 1 in [5]. Given this, we next show that $G_\varepsilon$ in (46) is in the support of $Q$. Fix $\delta > 0$ and denote by $B_\delta(c)$ the ball of radius $\delta$ centered in $c$. Define $W_l(G_\varepsilon) = \{G \in \mathbb{G} : G(z_l) - G(z_{l-1}) \in B_{\delta/(2N)}[G_\varepsilon(z_l) - G_\varepsilon(z_{l-1})]\}$ for $l = 1, \ldots, N$, with the convention that $G(z_0) := G(0^-) = 0$ if $z_1 = 0$ so that $G(z_1) - G(z_0) = G\{0\}$. Then $\bigcap_{l=1}^N W_l(G_\varepsilon) \subset \{G \in \mathbb{G} : \sup_{x \leq M} |G(x) - G_\varepsilon(x)| < \delta\}$. The sets $W_l(G_\varepsilon)$ are independent under $Q$ and each has positive probability. We conclude that, for any $\delta > 0$, $Q\{G \in \mathbb{G} : \sup_{x \leq M} |G(x) - G_\varepsilon(x)| < \delta\} \geq Q\{\bigcap_{l=1}^N W_l\} > 0$. The proof is then completed by taking $\varepsilon$ and $\delta$ such that $\varepsilon + \delta < \eta$. $\quad\square$

Now, relying on Theorem 2 and Lemma 11, we are in a position to provide the proofs of Theorems 4–7. Showing that for the specific kernels at issue (17) is met, represents a result of independent interest concerning small ball probabilities of mixtures with respect to CRMs; indeed, passing through Lemma 11, we actually show that (H1) is sufficient for (17), that is, for $\tilde{h}$ putting positive probability on uniform neighborhoods of $h_0$.

*Proof of Theorem 4.* The first step consists in verifying consistency with respect to hazards of mixture form. To this end we postulate the existence of a boundedly finite measure $\mu_0$ on $\mathbb{R}^+$ such that

$$(47) \qquad h_0(t) = \int_{\mathbb{R}^+} k(t, x)\mu_0(dx).$$

Clearly, $\mu_0$ has to be such that $\int_0^T h_0(t)\,dt \to +\infty$, as $T \to \infty$, in order to ensure the model to be properly defined. In the case of the DL kernel, (17) is a direct consequence of Lemma 11 since $h_0(t) = G_0(t)$ and $\tilde{h}(t) = \tilde{\mu}([0, t])$.

The consistency result clearly extends to all increasing hazard rates $h_0$ with $h_0(0) = 0$. To see this let $\mu_0$ be the measure associated to $h_0$. Then $\mu_0 \in \mathbb{M}$ since $\mu((0, \tau]) = h_0(\tau) \to 0$ as $\tau \to 0$ and $h_0(t) = \int \mathbb{I}_{(0 \leq x \leq t)}\mu_0(dx)$. Finally, note that the moment condition in (i) of Theorem 2 reduces to $\int_{\mathbb{R}^+} \mathbb{E}[\tilde{H}(t)]f_0(t)\,dt < \infty$ since, for any choice of $\lambda$ in (6) and for any large enough $t$, $\mathbb{E}[\tilde{H}(t)] > t$.

*Proof of Theorem 5.* As before, we first establish (17) for $h_0$ of mixture form (47) and assume $\tau$ to be fixed and (i) and (ii) to hold. Take $G \in \{G \in \mathbb{G} : \sup_{t \leq T+\tau} |G(x) - G_0(x)| < \delta\}$ and let $h_G$ be the corresponding hazard

rate. Then, one has

$$\sup_{t \leq T} |h_G(t) - h_0(t)|$$

$$= \sup_{t \leq T} \left| \int_{(t-\tau)^+}^{t+\tau} dG(x) - \int_{(t-\tau)^+}^{t+\tau} dG_0(x) \right|$$

$$= \sup_{t \leq T} |G(t+\tau) - G((t-\tau)^+) - G_0(t+\tau) + G_0((t-\tau)^+)|$$

$$\leq \sup_{t \leq T} |G(t+\tau) - G_0(t+\tau)| + \sup_{t \leq T} |G((t-\tau)^+) - G_0((t-\tau)^+)| \leq 2\delta.$$

Take $\delta$ such that $2\delta < \eta$, which yields $\Pi\{h : \sup_{0 < t \leq T} |h(t) - h_0(t)| < \eta\} \geq Q\{G \in \mathbb{G} : \sup_{x \leq T+\tau} |G(x) - G_0(x)| < \delta\}$, where we recall that $Q$ is the distribution induced on $\mathbb{G}$. The right-hand side has positive probability by Lemma 11 and, hence, (17) holds.

The next step consists in showing that any $h_0$, which is bounded Lipschitz continuous of constant $K > 0$ and satisfies (i) and (ii), can be approximated in the sup norm on $[0, T]$ to any order of accuracy by a rectangular mixture hazard (47) with a sufficiently small bandwidth $\tau$. To this end, define $h_m(t) = \int \mathbb{I}_{(|t-x| \leq \tau_m)} dG_m(x)$ with $\tau_m = m^{-\eta}$ ($m = 1, 2, \ldots$ and $\eta > 0$) and $dG_m(x) = \mathbb{I}_{(x < m)}(2\tau_m)^{-1} h_0(x) dx$. Note that $G_m \in \mathbb{G}$ for any integer $m$. Hence, we have

$$h_m(t) = \frac{1}{2\tau_m}(H_0(m \vee (t+\tau_m)) - H_0((t-\tau_m)^+))\mathbb{I}_{t < m+\tau_m}$$

and $h_m(t) \to h_0(t)$ for any $t$ as $m \to \infty$.

Next we apply the Arzelà–Ascoli theorem in order to obtain uniform convergence on a compact $[0, T]$. Hence we need to show that: (a) the sequence $\{h_m\}_{m \geq 1}$ is bounded on $[0, T]$ uniformly in $m$; (b) $\{h_m\}_{m \geq 1}$ is an equicontinuous sequence of functions on $[0, T]$. See Theorem 3 on page 270 in Feller [7]. Condition (a) is implied by $H_0$ being Lipschitz, which is guaranteed by the boundedness of $h_0$. Condition (b) boils down to showing that, to each $\varepsilon > 0$, there corresponds a $\delta > 0$ such that

$$(48) \qquad\qquad |t - s| < \delta \implies |h_m(t) - h_m(s)| < \varepsilon$$

for all large $m$. For simplicity we consider $\tau_m \leq s < t < m - \tau_m$. Then

$$|h_m(t) - h_m(s)| = \left| \frac{H_0(t+\tau_m) - H_0(t-\tau_m)}{2\tau_m} - \frac{H_0(s+\tau_m) - H_0(s-\tau_m)}{2\tau_m} \right|$$

$$= |h_0(t^*) - h_0(s^*)|$$

for some $t^*, s^*$ such that $t - \tau_m \leq t^* < t + \tau_m$ and $s - \tau_m \leq s^* < s + \tau_m$. Next, for $t^*, s^*$ running in these two intervals

$$|h_0(t^*) - h_0(s^*)| \leq \sup_{t^*, s^*} |h_0(t^*) - h_0(s^*)|$$

$$\leq K \sup_{t^*, s^*} |t^* - s^*| = K|t + \tau_m - (s - \tau_m)|$$

$$\leq K(\delta + 2\tau_m).$$

Finally, for given $\varepsilon$, choose $m_0$ large enough such that $T < m_0 - \tau_{m_0}$ and $\varepsilon/K - 2\tau_{m_0} > 0$. Then (48) is satisfied for $\delta < \varepsilon/K - 2\tau_{m_0}$, and (b) is proved.

Now fix $\eta > 0$. There exists $m$ such that $\sup_{0 < t \leq T} |h_m(t) - h_0(t)| < \eta/2$. For this $m$, take $G \in \{G \in \mathbb{G} : \sup_{x \leq T + \tau_m} |G(x) - G_m(x)| < \delta\}$ for $\delta < \eta/4$ and let $h_G$ be the corresponding hazard rate. Then, one has

$$\sup_{t \leq T} |h_G(t) - h_m(t)|$$

$$= \sup_{t \leq T} \left| \int_{(t-\tau)^+}^{t+\tau} dG(x) - \int_{(t-\tau)^+}^{t+\tau} dG_m(x) \right|$$

$$= \sup_{t \leq T} |G(t+\tau) - G((t-\tau)^+) - G_m(t+\tau) + G_0((t-\tau)^+)|$$

$$\leq \sup_{t \leq T} |G(t+\tau) - G_m(t+\tau)| + \sup_{t \leq T} |G((t-\tau)^+) - G_m((t-\tau)^+)| < \eta/2.$$

Such $m$ and $\delta$ yield $\Pi\{h : \sup_{0 < t \leq T} |h(t) - h_0(t)| < \eta\} \geq Q\{G \in \mathbb{G} : \sup_{x \leq T + \tau} |G(x) - G_0(x)| < \delta\} \times \pi\{\tau \in (0, \tau_m)\}$. The right-hand side has positive probability by Lemma 11 and the hypotheses on $\pi$. Hence, the proof is complete.

*Proof of Theorem 6.* As for the Ornstein–Uhlenbeck kernel, note that, since $\tilde{\mu}$ is a.s. discrete, $\tilde{h}$ is a shot-noise process with exponentially decaying shocks, that is, $\tilde{h}(t) = \sum_i J_i \sqrt{2\kappa} \exp\{-\kappa(t - X_i)\} \mathbb{I}_{(0 \leq X_i \leq t)}$, where the $J_i$'s and $X_i$'s are the random shocks and locations, respectively. We first aim at showing that any $h_0$ of the form (47) satisfying (i) and (ii) can be approximated in the sup norm on $[0, T]$ to any order of accuracy by a step-wise continuous function with a finite number of jumps. Let $G_\varepsilon$ be the step function defined in (46) and $h_\varepsilon$ the corresponding hazard rate $h_\varepsilon(t) = \sum_{l=1}^N j_l \sqrt{2\kappa} \exp\{-\kappa(t - z_l)\} \mathbb{I}_{(0 \leq z_l \leq t)}$. We first prove that

$$(49) \qquad \sup_{t \leq T} |h_0(t) - h_\varepsilon(t)| \leq \varepsilon\sqrt{2\kappa}$$

by determining a lower and an upper bound for the difference $h_0 - h_\varepsilon$. It turns out that the minimum distance $h_0 - h_\varepsilon$ is attained at one of the jump points $z_l$'s and a lower bound for $h_0(z_l)$ is obtained by moving the increment $G_0(z_i^-) - G_0(z_{i-1})$ near to the right of $z_{i-1}$ for any $i < l$. Setting $\Delta_i := \sqrt{2\kappa}[G_0(z_i^-) - G_0(z_{i-1})]$, we have

$$h_\varepsilon(z_l) - h_0(z_l) \leq h_\varepsilon(z_l) - \sum_{i=1}^l G_0\{z_i\}\sqrt{2\kappa}e^{-\kappa(z_l - z_i)} - \sum_{i=1}^l \Delta_i e^{-\kappa(z_l - z_{i-1})}$$

$$= \sum_{i=1}^{l} \Delta_i e^{-\kappa(z_l - z_i)} - \sum_{i=1}^{l} \Delta_i e^{-\kappa(z_l - z_{i-1})}$$

$$\leq \varepsilon\sqrt{2\kappa} \sum_{i=1}^{l} [e^{-\kappa(z_l - z_i)} - e^{-\kappa(z_l - z_{i-1})}] \leq \varepsilon\sqrt{2\kappa}.$$

As for the maximum of $h_0 - h_\varepsilon$, an upper bound for $h_0(t)$, with $t \in [z_l, z_{l+1})$, is obtained by moving the increment $G_0(z_i^-) - G_0(z_{i-1})$ near to the left of $z_i$ for $i \leq l$ and $G_0(t) - G_0(z_l)$ near to the left of $t$. Hence, we get

$$h_0(t) - h_\varepsilon(t) \leq h_\varepsilon(z_l) + [G_0(t) - G_0(z_l)]\sqrt{2\kappa} - h_\varepsilon(t) = [G_0(t) - G_0(z_l)]\sqrt{2\kappa}$$
$$\leq \varepsilon\sqrt{2\kappa}$$

and (49) is proved. Now, take $G$ such that $\sup_{x \leq T} |G_\varepsilon(x) - G(x)| < \delta$ and denote by $h_G(t) = \int_0^t \sqrt{2\kappa} \exp\{-\kappa(t-x)\} G(dx)$. We show that

$$(50) \qquad\qquad \sup_{t \leq T} |h_\varepsilon(t) - h_G(t)| \leq 2\delta\sqrt{2\kappa}.$$

Reasoning as for (49), the following bounds for $h_G(t)$ can be found

$$h_G(t) \leq \sqrt{2\kappa}(2\delta - \delta e^{-\kappa(t-z_1)^+}) + \sum_{i=1}^{N} j_i \sqrt{2\kappa} e^{-\kappa(t-z_i)} \mathbb{I}_{(0 \leq z_i \leq t)},$$

$$h_G(t) \geq \sqrt{2\kappa}(\delta e^{-\kappa t} - 2\delta e^{-\kappa(t-z_1)}) \mathbb{I}_{(z_1 \leq t)} + \sum_{i=1}^{N} j_i \sqrt{2\kappa} e^{-\kappa(t-z_i)} \mathbb{I}_{(0 \leq z_i \leq t)}$$

with $t \in [z_l, z_l + 1)$, $l = 0, \ldots, N-1$, $a^+ = a \vee 0$ and $\sum_{l=1}^{0} = 0$. Hence,

$$\sqrt{2\kappa}(\delta e^{-\kappa t} - 2\delta e^{-\kappa(t-z_1)}) \mathbb{I}_{(z_1 \leq t)} \leq h_\varepsilon(t) - h_G(t) \leq \sqrt{2\kappa}(2\delta - \delta e^{-\kappa(t-z_1)^+}),$$

which leads to the following bound in the sup norm

$$\sup_{t \leq T} |h_\varepsilon(t) - h_G(t)| \leq \max\{\sqrt{2\kappa}(2\delta - \delta e^{-\kappa(T-z_1)}), \sqrt{2\kappa}(2\delta - \delta e^{-\kappa T})\} \leq 2\delta\sqrt{2\kappa}.$$

Thus, (50) is proved. Now, by combining (49) and (50), for any $G \in \mathbb{G}$ such that $\sup_{x \leq T} |G(x) - G_\varepsilon(x)| < \delta$, we have

$$\sup_{t \leq T} |h_0(t) - h_G(t)| \leq \sup_{t \leq T} |h_0(t) - h_\varepsilon(t)| + \sup_{t \leq T} |h_\varepsilon(t) - h_G(t)| \leq (\varepsilon + 2\delta)\sqrt{2\kappa}.$$

Now, for any $\eta > 0$, take $\varepsilon$ and $\delta$ small enough such that $(\varepsilon + 2\delta)\sqrt{2\kappa} < \eta$. Hence, we obtain $\Pi\{h : \sup_{0 < t \leq T} |h(t) - h_0(t)| < \eta\} \geq Q\{G \in \mathbb{G} : \sup_{x \leq T} |G(x) - G_\varepsilon(x)| < \delta\}$. Note that the right-hand side has positive probability by Lemma 11 and, hence, (17) is proved. Now we show that any differentiable hazard rate such that $h_0(0) = 0$ and, according to condition (iii),

$$-h_0'(t)/h_0(t) \leq \kappa\sqrt{2\kappa},$$

can be represented as OU mixture (47) with absolutely continuous $\mu_0$. Define $u(x) = h_0'(x) + \kappa\sqrt{2\kappa}h_0(x)$. Such $u$ is a well defined Radon–Nikodým derivative of a boundedly finite measure with respect to the Lebesgue measure. This follows from the fact that $u(x) \geq 0$ for any $x$ by condition (iii) and that it is integrable in zero since $\int_0^\tau u(x)\,dx = h_0(\tau) + \kappa\sqrt{2\kappa}H_0(\tau)$. Set $dG_0(x) = u(x)\,dx$, where $G_0$ is the d.f. associated to $\mu_0$, and define $h_*(t) := \int_0^t \sqrt{2\kappa}e^{-\kappa(t-x)}\,dG_0(x)$. Then, both $h^*$ and $h_0$ are solution of the differential equation

$$dh(t) = -\kappa\sqrt{2\kappa}h(t)\,dt + dG_0(t) \qquad \text{with } h(0) = 0.$$

Thus, they coincide and the proof is complete.

*Proof of Theorem 7.* Assume first $h_0$ to be an exponential mixture (47) satisfying assumptions (i) and (ii). Then, $h_0$ is obviously strictly positive on $\mathbb{R}^+$. As for the uniform bound of $|\tilde{h}(t) - h_0(t)|$, it is useful to write $h_0(t) = \int_{\mathbb{R}^+} e^{-t/x}\,dG_0'(x)$ and $\tilde{h}(t) = \int_{\mathbb{R}^+} e^{-t/x}\,\tilde{\mu}'(dx)$, where

$$G_0'(x) = \int_0^x z^{-1}\,dG_0(z) \quad \text{and} \quad \tilde{\mu}'([0,x]) = \int_0^x z^{-1}\tilde{\mu}(dz).$$

Note that, by condition (ii) and the assumption that $\tilde{h}(0) < \infty$ a.s., $G_0'(x) < \infty$ and $\tilde{\mu}'([0,x]) < \infty$ a.s. for any finite $x$. Let $Q'$ denote the distribution induced on $\mathbb{G}$ by $\tilde{\mu}'$. One can check that, if Lemma 11 holds for $G_0$ and $\tilde{\mu}$, then, for any finite $M$ and $\eta > 0$,

$$(51) \qquad Q'\left\{ G' \in \mathbb{G} : \sup_{x \leq M} |G'(x) - G_0'(x)| < \eta \right\} > 0.$$

We now derive a bound for $|\tilde{h}(t) - h_0(t)|$ by exploiting the uniformly equicontinuity of the family of functions $\{e^{-t/x}, t \leq T\}$, as $x$ varies in the compact set $[0, M]$ for any $T < \infty$ and $M < \infty$. In fact, given $\gamma > 0$, the Arzelà–Ascoli theorem ensures the existence of finitely many points $t_1, \ldots, t_m$ such that, for any $t \leq T$, there is an index $i$ for which

$$(52) \qquad \sup_{x \leq M} |e^{-t/x} - e^{-t_i/x}| \leq \gamma.$$

Now, note that condition (ii) and the assumption that $\tilde{h}(0) < \infty$ a.s., imply that: (i) for any $\varepsilon_1 > 0$, there exists $M_1 < \infty$ large enough such that $\int_{M_1}^\infty dG_0'(x) < \varepsilon_1$; (ii) for any $\varepsilon_2 > 0$, $\exists M_2 < \infty$ large enough such that $Q'\{G' : \int_{M_2}^\infty dG'(x) < \varepsilon_2\} > 0$. At this point, take $M = M_1 \vee M_2$ and note that $\{G' : \int_{M_2}^\infty dG'(x) < \varepsilon_2\} \subseteq A(M, \varepsilon_2)$, where $A(M, \varepsilon_2) := \{G' : \int_M^\infty dG'(x) < \varepsilon_2\}$. Finally, define

$$B(M, \varepsilon_3) := \left\{ G' \in \mathbb{G} : \sup_{x \leq M} |G'(x) - G_0'(x)| < \varepsilon_3 \right\},$$

which is a set of positive probability for any $\varepsilon_3 > 0$ by (51), and note that $Q'[A(M, \varepsilon_2) \cap B(M, \varepsilon_3)] > 0$ by the independence of $\tilde{\mu}((0, M])$ and $\tilde{\mu}([M, \infty))$. Take now $G' \in A(M, \varepsilon_2) \cap B(M, \varepsilon_3)$ and let $h_{G'}(t) = \int_{\mathbb{R}^+} e^{-t/x} dG'(x)$. Then, for an arbitrary $t \leq T$, choose the appropriate $t_i$ such that (52) holds and write

$$|h_{G'}(t) - h_0(t)| \leq \int_{\mathbb{R}^+} |e^{-t/x} - e^{-t_i/x}| \, dG'(x) + \int_{\mathbb{R}^+} |e^{-t/x} - e^{-t_i/x}| \, dG'_0(x)$$

$$+ \left| \int_{\mathbb{R}^+} e^{-t_i/x} \, dG'(x) - \int_{\mathbb{R}^+} e^{-t_i/x} \, dG'_0(x) \right| := I_1 + I_2 + I_3.$$

As for $I_1$, we have

$$I_1 \leq \gamma G'(M) + \int_M^\infty |e^{-t/x} - e^{-t_i/x}| \, dG'(x)$$

$$\leq \gamma G'(M) + 2 \int_M^\infty dG'(x) \leq \gamma[G'_0(M) + \varepsilon_3] + 2\varepsilon_2,$$

where we used the fact that $h_0$ and $\tilde{h}$ are decreasing in the second step. Similar arguments lead to $I_2 \leq \gamma G'_0(M) + 2\varepsilon_1$. Concerning $I_3$, write

$$I_3 \leq \left| \int_0^M e^{-t_i/x} \, dG'(x) - \int_0^M e^{-t_i/x} \, dG'_0(x) \right|$$

$$+ \left| \int_M^\infty e^{-t_i/x} \, dG'(x) - \int_M^\infty e^{-t_i/x} \, dG'_0(x) \right|$$

$$\leq \left| \int_0^M e^{-t_i/x} G'(dx) - \int_0^M e^{-t_i/x} G'_0(dx) \right| + \varepsilon_2 + \varepsilon_1 \leq \varepsilon_3 + \varepsilon_2 + \varepsilon_1,$$

where, in the last step, we have exploited the fact that $G'$ belongs also to a weak neighborhood of $G'_0$ of radius $\varepsilon_3$, when one reasons in terms of finite measures over $[0, M]$. Summing up, we have obtained $|h_{G'}(t) - h_0(t)| \leq 2\gamma G'_0(M) + \gamma \varepsilon_3 + \varepsilon_3 + 3(\varepsilon_2 + \varepsilon_1)$, where $G'_0(M)$ is a finite constant. Hence, we are able to state that, for a given $\eta$, it is always possible to choose $\gamma, \varepsilon_1, \varepsilon_2$ and $\varepsilon_3$ such that $|h_{G'}(t) - h_0(t)| \leq \eta$ for $G'$ in a set of positive probability. To see this, set $\varepsilon_1$ and $\varepsilon_2$ such that $3(\varepsilon_1 + \varepsilon_2) < \eta/4$, then determine $M = M_1 \vee M_2$; since, for such $M$, we have $G'_0(M) < \infty$, set $\gamma$ such that $2\gamma G'_0(M) < \eta/4$; for such $\gamma$ set $\varepsilon_3$ such that $\varepsilon_3(\gamma + 1) < \eta/4$.

The next step consists in establishing that any function completely monotone function $\varphi$ on $\mathbb{R}^+$ such that $\varphi(0) < \infty$ is of the form

$$\varphi(t) = \int_0^\infty x^{-1} e^{-t/x} dG(x), \tag{53}$$

where $G \in \mathbb{G}$, that is, it is an exponential mixture with respect to a boundedly finite measure. The starting point is the fundamental result of Bernstein, which characterizes completely monotone functions as mixtures, the mixing measures being probability measures. For our needs it is more convenient to resort to the version of Bernstein's result as formulated in Theorem 1(a) on page 439 in Feller [7]: a function $\varphi$ on $(0, \infty)$ is completely monotone if and only if it is of the form

$$(54) \qquad \varphi(\lambda) = \int_0^\infty e^{-\lambda x} U(dx),$$

where $U \in \mathbb{M}$. Without loss of generality, we may assume that $\varphi \sim L(1/\tau)$, where $L$ is a slowly varying function at infinity. This clearly covers the case of $\varphi$'s such that $\varphi(0) < \infty$, which is one of our assumptions. At this point we resort to a suitable Tauberian theorem (see Theorem on page 445 in [7]), which allows to deduce the behavior at infinity of $U$ in (54) from the behavior in zero of $\varphi$. Hence, we have $U(t) \sim L(t)$, as $t \to \infty$. Let now $T(x) = 1/x$ and denote with $U \circ T^{-1}$ the image measure of $U$ by $T$. We can write

$$\varphi(\lambda) = \int_0^\infty e^{-\lambda/y}(U \circ T^{-1})(dy) = \int_0^\infty y^{-1} e^{-\lambda/y} y(U \circ T^{-1})(dy)$$

and define $G(dy) \equiv y(U \circ T^{-1})(dy)$. For simplicity, we assume that $U(x)$ has an ultimately monotone derivative, that is, $U(dx) = u(x)\, dx$ with $u(x)$ monotone in some interval $(x_0, \infty)$. Then

$$G(\tau) = \int_0^\tau y(U \circ T^{-1})(dy) = \int_{1/\tau}^\infty x^{-1} u(x)\, dx$$

for sufficiently small $\tau$. We aim at showing that $G(\tau) \to 0$ as $\tau \to 0$. In fact, $U(t) \sim L(t)$ implies that, for any $\varepsilon > 0$, $u(t) = o(t^{\varepsilon-1} L(t))$ as $t \to \infty$. Otherwise, if $u(t) \sim K t^{\varepsilon^*-1} L(t)$ for some $\varepsilon^* > 0$ and constant $K$, then $U(t) \sim (K/\varepsilon^*) t^{\varepsilon^*} L(t)$ (see the lemma after Theorem 4 on page 446 in [7]) which, in turn, contradicts $U(t) \sim L(t)$. Next we have

$$\int_{1/\tau}^\infty x^{-1} \frac{u(x)}{\tau^{1-\varepsilon} L(1/\tau)}\, dx = \int_1^\infty y^{-1} \frac{u(y/\tau)}{\tau^{1-\varepsilon} L(1/\tau)}\, dy \to 0 \qquad \text{as } \tau \to 0,$$

where the integrand is monotone and it remains bounded as $\tau \to 0$. Thus, $G(\tau) = o(\tau^{1-\varepsilon} L(1/\tau))$ for any $\varepsilon > 0$, and in particular for $\varepsilon < 1$, from which the desired result follows. We have then established that any completely monotone function $\varphi$ such that $\varphi(0) < \infty$ is of the form (53).

Finally, the fact that the moment condition (ii) in Theorem 2 reduces to $\int t f_0(t)\, dt < \infty$ follows from the fact that the function $t \mapsto \tilde{h}(t)$ is a.s. decreasing. Hence, the proof is complete.

**Further results and proofs for Section 4.**

*Compensated Poisson random measures.* In order to prove the results concerning functionals of hazard rates, we will often work with the *compensated Poisson random measure* canonically associated to a Poisson measure $\tilde{N}$ with intensity $\nu$. This object is written $\tilde{N}^c = \{\tilde{N}^c(A) : A \in \mathscr{B}(\mathbb{R}^+) \otimes \mathscr{X}\}$ and is defined as the unique CRM on $(\mathbb{R}^+ \times \mathbb{X}, \mathscr{B}(\mathbb{R}^+) \otimes \mathscr{X})$ such that

$$\tilde{N}^c(A) = \tilde{N}(A) - \nu(A)$$

for every set $A$ of finite $\nu$-measure. For every $g \in L^2(\nu)$, we denote by

$$\tilde{N}^c(g) = \int_{\mathbb{R}^+ \times \mathbb{X}} g(s, x) \tilde{N}^c(ds, dx)$$

the Wiener–Itô integral of $g$ with respect to $\tilde{N}^c$. Observe that, for every $g \in L^2(\nu)$, $\tilde{N}^c(g)$ is a centered and square integrable random variable with an infinitely divisible law. In particular, for every $\lambda \in \mathbb{R}$,

$$(55) \qquad \mathbb{E}[e^{i\lambda \tilde{N}^c(g)}] = \exp\left\{ \int_{\mathbb{R}^+ \times \mathbb{X}} [e^{i\lambda g(s,x)} - 1 - i\lambda g(s,x)] \nu(ds, dx) \right\}.$$

Also for every $f, g \in L^2(\nu)$, one has the fundamental *isometric property*

$$(56) \qquad \mathbb{E}[\tilde{N}^c(f) \tilde{N}^c(g)] = \int_{\mathbb{R}^+ \times \mathbb{X}} f(s, x) g(s, x) \nu(ds, dx) := (f, g)_{L^2(\nu)}.$$

Note that (35), (55) and (56) imply that, for every $g \in L^2(\nu) \cap L^1(\nu)$,

$$\mathbb{E}[\tilde{N}(g)] = \int_{\mathbb{R}^+ \times \mathbb{X}} g(s, x) \nu(ds, dx)$$

$$\mathrm{Var}[\tilde{N}(g)] = \mathrm{Var}[\tilde{N}^c(g)] = \int_{\mathbb{R}^+ \times \mathbb{X}} g(s, x)^2 \nu(ds, dx).$$

*Limit theorems for shifted measures.* In this section we prove a series of preliminary CLTs, involving random hazard rates that are obtained from $\tilde{h}^{n,*}$ [as defined in (18)] by adding fixed atoms to the underlying CRM $\tilde{\mu}^{n,*}$. The notation and framework are those of Sections 2 and 4.1.

Fix a natural number $k \geq 1$, along with points $x_1, \ldots, x_k \in \mathbb{X}$ such that $x_i \neq x_j$ for every $i \neq j$, and positive coefficients $z_1, \ldots, z_k \in \mathbb{R}^+$. We define the discrete measure $\Delta(\cdot)$, on $(\mathbb{X}, \mathscr{X})$ as follows:

$$(57) \qquad \Delta(B) = \sum_{j=1}^{k} z_j \delta_{x_j}(B), \qquad B \in \mathscr{X},$$

where $\delta_y$ stands for the Dirac mass concentrated at $y$. Now set $\tilde{\mu}^{n,*}_{\Delta}(B) = \tilde{\mu}^{n,*}(B) + \Delta(B)$, for $B \in \mathscr{X}$, where $\tilde{\mu}^{n,*}$ is the CRM appearing in (12), and

also

$$\tilde{h}_\Delta^{n,*}(t) = \int_{\mathbb{X}} k(t,x) \tilde{\mu}_\Delta^{n,*}(dx) = \tilde{h}^{n,*}(t) + \sum_{j=1}^{k} z_j k(t,x_j)$$

(58)

$$\tilde{H}_\Delta^{n,*}(T) = \int_0^T \tilde{h}_\Delta^{n,*}(t)\,dt = \tilde{H}^{n,*}(T) + \sum_{j=1}^{T} z_j \int_0^T k(t,x_j)\,dt.$$

Note that, with the notation introduced in (58), one has that the cumulative hazard rate with fixed random jumps $\tilde{H}_{\Delta^{n,*}}(T)$ is indeed such that $\tilde{H}_{\Delta^{n,*}}(T) = \tilde{H}_{\Delta^{n,*}}^{n,*}(T)$; however, this heavy notation is avoided henceforth.

Our aim is now to establish CLTs for linear and quadratic functionals of the transformed hazard rate $\tilde{h}_\Delta(\cdot)$. These results represent the "deterministic skeleton" upon which the conditioned CLTs of Section 4 are constructed. Note that the random measure $\tilde{\mu}_\Delta^{n,*}$ is a CRM with fixed atoms (given by the points $x_1, \ldots, x_k$), so that one cannot apply directly the theories developed in [27, 28]. An integer $n \geq 0$ is fixed for the rest of the section.

PROPOSITION 12. *Suppose that points* (i) *and* (ii) *in the statement of Theorem 8 are satisfied for $n \geq 0$. Assume moreover that*

(59) $$\lim_{T \to +\infty} C_0(n,k,T) \times \sum_{j=1}^{k} z_j \int_0^T k(t,x_j)\,dt = m(n,\Delta,k) \in [0,+\infty).$$

*Then, letting $X \sim \mathcal{N}(m(n,\Delta,k), \sigma_0^2(n,k))$, we have*

$$C_0(n,k,T) \times [\tilde{H}_\Delta^{n,*}(T) - \mathbb{E}[\tilde{H}^{n,*}(t)^2]] \xrightarrow{\text{law}} X.$$

Before proving the result, it is worth pointing out that (59) only involves deterministic quantities, and also that we do not suppose (29) to hold.

PROOF. First, write

$$C_0(n,k,T) \times [\tilde{H}_\Delta^{n,*}(T) - \mathbb{E}[\tilde{H}^{n,*}(t)^2]]$$

$$= C_0(n,k,T) \times [\tilde{H}^{n,*}(T) - \mathbb{E}[\tilde{H}^{n,*}(t)^2]]$$

$$+ C_0(n,k,T) \times \sum_{j=1}^{k} z_j \int_0^T k(t,x_j)\,dt.$$

Now observe that $\tilde{H}^{n,*}(T)$ is the cumulative hazard rate obtained from a CRM with intensity $\nu^{n,*}$. As a consequence, according to Theorem 1 in [27], whenever conditions (i) and (ii) of Theorem 8 are verified, one has that the sequence $C_0(n,k,T) \times [\tilde{H}^{n,*}(T) - \mathbb{E}[\tilde{H}^{n,*}(t)^2]]$ converges in law to a Gaussian random variable with variance $\sigma_0^2(n,k)$. Since (59) holds by assumption, and since $m(n,\Delta,k)$ is deterministic, the conclusion follows.  $\square$

Now, define the kernel $k_{T,\Delta}^{(4)}$ as in (26), by simply replacing $\Delta^{n,*}$ with $\Delta$, that is:

$$k_{T,\Delta}^{(4)}(s,x) := \sum_{j=1}^{k} k_T^{(1)}(s,x;z_j,x_j), \qquad (s,x) \in \mathbb{R}^+ \times \mathbb{X}.$$

PROPOSITION 13. *Suppose that all the assumptions in the statement of Theorem 9 are satisfied, except for points 5., 6. and 7., which are replaced, respectively, by*

5b.    $C_1^2(n,k,T)\|k_T^{(2)} + 2k_{T,n}^{(3)} + 2k_{T,\Delta}^{(4)}\|_{L^2(\nu^{n,*})}^2 \to \sigma_4^2(n,\Delta,k) \geq 0;$

6b.    $C_1^3(n,k,T)\|k_T^{(2)} + 2k_{T,n}^{(3)} + 2k_{T,\Delta}^{(4)}\|_{L^3(\nu^{n,*})}^3 \to 0;$

7b.    $\displaystyle\lim_{T \to +\infty} \frac{C_1(n,k,T)}{T} \int_0^T \left(\sum_{j=1}^{k} z_j k(t,x_j)\right)^2 dt = v(n,\Delta,k) \in [0,+\infty).$

*Then, letting $X \sim \mathscr{N}(v(n,\Delta,k), \sigma_1^2(n,k) + \sigma_4^2(n,\Delta,k))$, we have*

$$
\begin{aligned}
(60) \quad C_1(n,k,T) \times \Bigg\{ &\frac{1}{T}\int_0^T \tilde{h}_\Delta^{n,*}(t)^2\, dt \\
&- \sum_{j=1}^{k} \frac{2z_j}{T}\int_0^T \mathbb{E}[\tilde{h}^{n,*}(t)]k(t,x_j)\, dt \\
&- \frac{1}{T}\int_0^T \mathbb{E}[\tilde{h}^{n,*}(t)^2]\, dt \Bigg\} \xrightarrow{\text{law}} X.
\end{aligned}
$$

PROOF.    Denote by $\tilde{N}^{n,*}$ the Poisson measure on $\mathbb{R}^+ \times \mathbb{X}$, with intensity $\nu^{n,*}$, determining $\tilde{\mu}^{n,*}$. We also write $\tilde{N}^{c;n,*}$ to indicate the compensated Poisson measure associated with $\tilde{N}^{n,*}$. First observe that, by, for example, Lemma 1 in [27],

$$\frac{1}{T}\int_0^T \tilde{h}^{n,*}(t)\sum_{j=1}^{k} z_j k(t,y_j)\, dt$$

$$\stackrel{\text{law}}{=} \frac{1}{T}\int_0^T \tilde{N}^{c;n,*}((\cdot)k(t,\cdot))\sum_{j=1}^{k} z_j k(t,x_j)\, dt$$

$$+ \frac{1}{T}\int_0^T \int_{\mathbb{R}^+ \times \mathbb{X}} s k(t,x)\nu^{n,*}(ds,dx)\sum_{j=1}^{k} z_j k(t,x_j)\, dt$$

$$= \tilde{N}^{c;n,*}(k_{T,\Delta}^{(4)}) + \sum_{j=1}^{k} \frac{z_j}{T} \int_0^T \mathbb{E}[\tilde{h}^{n,*}(t)]k(t,x_j)\,dt,$$

where $\tilde{N}^{c;n,*}((\cdot)k(t,\cdot)) := \int_{\mathbb{R}^+\times\mathbb{X}} sk(t,x)\tilde{N}^{c;n,*}(ds,dx)$. From the previous relations, we deduce

$$
\frac{1}{T}\int_0^T \tilde{h}_\Delta^{n,*}(t)^2\,dt \stackrel{\text{law}}{=} \frac{1}{T}\int_0^T (\tilde{h}^{n,*}(t))^2\,dt + \frac{1}{T}\int_0^T \left(\sum_{j=1}^k z_j k(t,x_j)\right)^2\,dt
$$
(61)
$$
+ \tilde{N}^{c;n,*}(2k_{T,\Delta}^{(4)}) + 2\sum_{j=1}^k \frac{z_j}{T}\int_0^T \mathbb{E}[\tilde{h}^{n,*}(t)]k(t,x_j)\,dt.
$$

Now apply the calculations contained in [27], Section 5.3, to deduce that

$$
\frac{1}{T}\int_0^T \tilde{h}^{n,*}(t)^2\,dt - \frac{1}{T}\int_0^T \mathbb{E}[\tilde{h}^{n,*}(t)^2]\,dt
$$
(62)
$$
\stackrel{\text{law}}{=} \tilde{N}^{c;n,*}(k_T^{(2)} + 2k_{T,n}^{(3)}) + I_2(k_T^{(1)}),
$$

where $I_2$ stands for a double Poisson integral with respect to $\tilde{N}^{c;n,*}$ (see [27] or [28] for further details). From the last formula and from (61) we infer that the expression in (60) has indeed the same law as

$$
C_1(n,k,T)[\tilde{N}^{c;n,*}(k_T^{(2)} + 2k_{T,n}^{(3)} + 2k_{T,\Delta}^{(4)}) + I_2(k_T^{(1)})]
$$
$$
+ C_1(n,k,T)\frac{1}{T}\int_0^T \left(\sum_{j=1}^k z_j k(t,x_j)\right)^2\,dt.
$$

To justify the operation of "plugging" the equality in law (62) into (61), one can use the more general relation: $\tilde{h}^{n,*}(t) \stackrel{\text{law}}{=} \int_{\mathbb{R}^+\times\mathbb{X}} sk(t,x)\tilde{N}^{n,*}(ds,dx)$, where the equality holds in the sense of stochastic processes; see again Lemma 1 in [27]. Now we can apply directly Theorem 3 in [28] to deduce that, under the assumptions in the statement, the pair

$$
C_1(n,k,T)(\tilde{N}^{c;n,*}(k_T^{(2)} + 2k_{T,n}^{(3)} + 2k_{\Delta,T}^{(4)}), I_2(k_T^{(1)}))
$$

converges in law to $(N,N')$ where $N,N'$ are two independent centered Gaussian random variables with variances given, respectively, by $\sigma_1^2(n,k)$ and $\sigma_4^2(n,\Delta,k)$. Since $v(n,\Delta,k)$ is deterministic, the conclusion follows. $\square$

PROPOSITION 14. *Suppose that the assumptions of Propositions 12 and 13 are verified. Assume also that points 1. and 2. in the statement of Theorem 10 hold and that point 3. in the same statement is replaced by*

3b: $\|C_1(n,k,T)(k_T^{(2)} + 2k_{T,n}^{(3)} + 2k_{T,\Delta}^{(4)}) - \delta(n,k)C_0(n,k,T)k_T^{(0)}\|_{L^2(\nu^{n,*})}^2$

$\qquad \rightarrow \sigma_5^2(n,\Delta,k) \geq 0.$

*Then, letting* $X \sim \mathcal{N}(v(n, \Delta, k) - \delta(n, k)m(n, \Delta, k), \sigma_1^2(n, k) + \sigma_5^2(n, \Delta, k))$,

$$
\begin{aligned}
V(T) := C_1(n, k, T) \Bigg\{ &\frac{1}{T} \int_0^T \left[ \tilde{h}_\Delta^{n,*}(t) - \frac{\tilde{H}_\Delta^{n,*}(T)}{T} \right]^2 dt \\
&- \sum_{j=1}^k \frac{2z_j}{T} \int_0^T \mathbb{E}[\tilde{h}^{n,*}(t)]k(t, y_j)\, dt \\
&- \frac{1}{T} \int_0^T \mathbb{E}[\tilde{h}^{n,*}(t)^2]\, dt + \frac{\mathbb{E}[\tilde{H}^{n,*}(T)]^2}{T^2} \Bigg\} \xrightarrow{\text{law}} X.
\end{aligned}
$$

(63)

PROOF.  Throughout the proof, we use the symbol $A(T) \overset{\mathbb{P}}{\approx} B(T)$ to indicate that $A(T) - B(T)$ converges to zero in probability. First observe that

$$
\begin{aligned}
&C_1(n, k, T) \left( \frac{\tilde{H}_\Delta^{n,*}(T)}{T} \right)^2 \\
&= \frac{C_1(n, k, T)\{C_0(n, k, T)[\tilde{H}_\Delta^{n,*}(T) - \mathbb{E}(\tilde{H}^{n,*}(T))]\}^2}{T^2 C_0(n, k, T)^2} \\
&\quad + \frac{C_1(n, k, T)}{T^2} \mathbb{E}(\tilde{H}^{n,*}(T))^2 \\
&\quad + 2\frac{C_1(n, k, T)}{T^2} \mathbb{E}(\tilde{H}^{n,*}(T))[\tilde{H}_\Delta^{n,*}(T) - \mathbb{E}(\tilde{H}(T))].
\end{aligned}
$$

Since point 1 in the statement of Theorem 10 is verified, and since the assumptions of Proposition 12 are in order, we deduce

$$
\frac{C_1(k, T)}{T^2 C_0(k, T)^2} \{C_0(k, T)[\tilde{H}(T) - \mathbb{E}(\tilde{H}(T))]\}^2 \overset{\mathbb{P}}{\to} 0.
$$

Moreover, point 2 in the statement of Theorem 10 yields that, as $T \to +\infty$,

$$
\begin{aligned}
&\frac{2C_1(n, k, T)}{T^2} \mathbb{E}(\tilde{H}^{n,*}(T))[\tilde{H}_\Delta^{n,*}(T) - \mathbb{E}(\tilde{H}^{n,*}(T))] \\
&\overset{\mathbb{P}}{\approx} \delta(n, k)C_0(n, k, T)[\tilde{H}_\Delta^{n,*}(T) - \mathbb{E}(\tilde{H}(T))].
\end{aligned}
$$

Now consider the functional $V(T)$ defined in (63). By reasoning as in the proofs of Propositions 12 and 13 we deduce that

$$
\begin{aligned}
V(T) \overset{\mathbb{P}}{\approx}\ &C_1(n, k, T) \\
&\times \Bigg\{ \frac{1}{T} \int_0^T [\tilde{h}_\Delta^{n,*}(t)^2]\, dt - \delta(n, k)C_0(n, k, T)[\tilde{H}_\Delta^{n,*}(T) - \mathbb{E}(\tilde{H}(T))] \\
&\quad - \sum_{j=1}^k \frac{z_j}{T} \int_0^T \mathbb{E}[\tilde{h}^{n,*}(t)]k(t, y_j)\, dt - \frac{1}{T} \int_0^T \mathbb{E}[\tilde{h}^{n,*}(t)^2]\, dt \Bigg\}
\end{aligned}
$$

$$\overset{\text{law}}{=} \tilde{N}^{c;n,*}(C_1(n,k,T)(k_T^{(2)} + 2k_{T,n}^{(3)} + 2k_{T,\Delta}^{(4)}) - \delta(n,k)C_0(n,k,T)k_T^{(0)})$$

$$+ I_2(C_1(n,k,T)k_T^{(1)}) - \delta(n,k)C_0(n,k,T)\sum_{j=1}^{k} z_j \int_0^T k(t,x_j)\,dt$$

$$+ \frac{C_1(n,k,T)}{T} \int_0^T \left(\sum_{j=1}^{k} z_j k(t,x_j)\right)^2 dt$$

and the conclusion is again obtained from Theorem 3 in [28]. □

*Proofs of Theorems 8, 9 and 10.* For the sake of brevity, we only provide the complete proof of Theorem 8. From point (i) of Theorem 1 we deduce

$$\mathbb{E}[e^{i\lambda C_0(n,k,T)[\tilde{H}(T)-\mathbb{E}[\tilde{H}^{n,*}(T)]]}|\mathbf{X}^* = (x_1,\ldots,x_k), \mathbf{Y}]$$

$$= \mathbb{E}[\exp(i\lambda C_0(n,k,T)[\tilde{H}_{\Delta\mathbf{J}}^{n,*}(T) - \mathbb{E}[\tilde{H}^{n,*}(T)]])],$$

where

$$\tilde{H}_{\Delta\mathbf{J}}^{n,*}(T) := (\tilde{H}_{\Delta^{n,*}}(T)|\mathbf{X}^* = (x_1,\ldots,x_k)) = \tilde{H}^{n,*}(T) + \sum_{i=1}^{k} J_i \int_0^T k(t,x_i)\,dt.$$

$\tilde{H}^{n,*}(T)$ and $\tilde{H}_{\Delta^{n,*}}(T)$ are defined in (19), and the jump vector $\mathbf{J} = (J_1,\ldots,J_k)$ is independent of $\tilde{H}^{n,*}(T)$ and with law given by (14). The previous relations and independence yield that

(64)
$$\mathbb{E}[e^{i\lambda C_0(n,k,T)\{\tilde{H}(T)-\mathbb{E}[\tilde{H}^{n,*}(T)]\}}|\mathbf{J} = (z_1,\ldots,z_k), \mathbf{X}^* = (x_1,\ldots,x_k), \mathbf{Y}]$$

$$= \mathbb{E}[e^{i\lambda C_0(n,k,T)\{\tilde{H}_\Delta^{n,*}(T)-\mathbb{E}[\tilde{H}^{n,*}(T)]\}}],$$

where $\tilde{H}_\Delta^{n,*}(T)$ is defined in (58) and $\Delta$ is given by (57). Now suppose that the assumptions of Theorem 8 are met [in particular, (29)]. This implies that there exists a set $\Omega'$ of $\mathbb{P}$-probability one such that, for every $\omega \in \Omega'$, the probability $B \mapsto \mathbb{P}\{\mathbf{X}^* \in B|\mathbf{Y}\}$ has support contained in the set of those vectors $(x_1,\ldots,x_k)$ such that, for every fixed $(z_1,\ldots,z_k)$ in the support of the law of $\mathbf{J}$, the cumulative hazard rate $\tilde{H}_\Delta^{n,*}(T)$ appearing in (64) verifies the assumptions of Proposition 12 [in particular, (59) holds with $m(n,\Delta,k) = m(n,\Delta^{n,*},k)$]. This yields, for all such $(x_1,\ldots,x_k)$ and $(z_1,\ldots,z_k)$,

$$\mathbb{E}[e^{i\lambda C_0(n,k,T)[\tilde{H}(T)-\mathbb{E}[\tilde{H}^{n,*}(T)]]}|\mathbf{J} = (z_1,\ldots,z_k), \mathbf{X}^* = (x_1,\ldots,x_k), \mathbf{Y}]$$

$$\underset{T\to\infty}{\longrightarrow} \exp\left(i\lambda m(n,\Delta^{n,*},k) - \frac{\lambda^2}{2}\sigma_0^2(n,k)\right).$$

To conclude, it is sufficient to use the Dominated Convergence Theorem for conditional expectations, to obtain that, a.s.-$\mathbb{P}$,

$$\mathbb{E}[e^{i\lambda C_0(n,k,T)[\tilde{H}(T)-\mathbb{E}[\tilde{H}^{n,*}(T)]]}|\mathbf{Y}]$$

$$= \mathbb{E}[\mathbb{E}[e^{i\lambda C_0(n,k,T)[\tilde{H}(T) - \mathbb{E}[\tilde{H}^{n,*}(T)]]}|\mathbf{X}^*, \mathbf{Y}]|\mathbf{Y}]$$

$$= \mathbb{E}[\mathbb{E}[\mathbb{E}[e^{i\lambda C_0(n,k,T)[\tilde{H}(T) - \mathbb{E}[\tilde{H}^{n,*}(T)]]}|\mathbf{J}, \mathbf{X}^*, \mathbf{Y}]|\mathbf{X}^*, \mathbf{Y}]|\mathbf{Y}]$$

$$\xrightarrow[T\to\infty]{} \mathbb{E}\Big[\mathbb{E}\Big[\exp\Big(i\lambda m(n, \Delta^{n,*}, k)(\omega) - \frac{\lambda^2}{2}\sigma_0^2(n,k)\Big)\Big|\mathbf{X}^*, \mathbf{Y}\Big]\Big|\mathbf{Y}\Big]$$

$$= \mathbb{E}\Big[\exp\Big(i\lambda m(n, \Delta^{n,*}, k)(\omega) - \frac{\lambda^2}{2}\sigma_0^2(n,k)\Big)\Big|\mathbf{Y}\Big],$$

thus completing the proof.

The proofs of Theorem 9 and 10 can be obtained by using exactly the same line of reasoning and by applying, respectively, Propositions 13 and 14.

**Acknowledgments.** The authors are grateful to an Associate Editor and two referees for important comments and suggestions, which led to a substantial improvement of the paper. Moreover, S. Ghosal and S. G. Walker are gratefully acknowledged for some useful discussions.

## REFERENCES

[1] BARRON, A., SCHERVISH, M. J. and WASSERMAN, L. (1999). The consistency of distributions in nonparametric problems. *Ann. Statist.* **27** 536–561. MR1714718

[2] BRIX, A. (1999). Generalized gamma measures and shot-noise Cox processes. *Adv. in Appl. Prob.* **31** 929–953. MR1747450

[3] DALEY, D. and VERE-JONES, D. J. (1988). *An Introduction to the Theory of Point Processes*. Springer, New York. MR0950166

[4] DOKSUM, K. (1974). Tailfree and neutral random probabilities and their posterior distributions. *Ann. Probab.* **2** 183–201. MR0373081

[5] DRĂGHICI, L. and RAMAMOORTHI, R. V. (2003). Consistency of Dykstra-Laud priors. *Sankhyā* **65** 464–481. MR2028910

[6] DYKSTRA, R. L. and LAUD, P. (1981). A Bayesian nonparametric approach to reliability. *Ann. Statist.* **9** 356–367. MR0606619

[7] FELLER, W. (1971). *An Introduction to Probability Theory and Its Applications*. Vol. II, 3rd ed. Wiley, New York. MR0270403

[8] GHOSAL, S., GHOSH, J. K. and RAMAMOORTHI, R. V. (1999). Posterior consistency of Dirichlet mixtures in density estimation. *Ann. Statist.* **27** 143–158. MR1701105

[9] GHOSH, J. K. and RAMAMOORTHI, R. V. (1995). Consistency of Bayesian inference for survival analysis with or without censoring. In *Analysis of Censored Data (Pune, 1994/1995). IMS Lecture Notes Monogr. Ser.* **27** 95–103. IMS, Hayward, CA. MR1483342

[10] GHOSH, J. K. and RAMAMOORTHI, R. V. (2003). *Bayesian Nonparametrics*. Springer, New York. MR1992245

[11] HJORT, N. L. (1990). Nonparametric Bayes estimators based on beta processes in models for life history data. *Ann. Statist.* **18** 1259–1294. MR1062708

[12] HO, M.-W. (2006). A Bayes method for a monotone hazard rate via *S*-paths. *Ann. Statist.* **34** 820–836. MR2283394

[13] ISHWARAN, H. and JAMES, L. F. (2004). Computational methods for multiplicative intensity models using weighted gamma processes: Proportional hazards, marked point processes, and panel count data. *J. Amer. Statist. Assoc.* **99** 175–190. MR2054297

[14] JAMES, L. F. (2003). Bayesian calculus for gamma processes with applications to semiparametric intensity models. *Sankhyā* **65** 179–206. MR2016784

[15] JAMES, L. F. (2005). Bayesian Poisson process partition calculus with an application to Bayesian Lévy moving averages. *Ann. Statist.* **33** 1771–1799. MR2166562

[16] KABANOV, Y. (1975). On extendend stochastic integrals. *Teor. Verojatnost. i Primenen.* **20** 725–737. MR0397807

[17] KALLENBERG, O. (1986). *Random Measures*, 4th ed. Akademie-Verlag, Berlin. MR0854102

[18] KIM, Y. (1999). Nonparametric Bayesian estimators for counting processes. *Ann. Statist.* **27** 562–588. MR1714717

[19] KIM, Y. (2003). On the posterior consistency of mixtures of Dirichlet process priors with censored data. *Scand. J. Statist.* **30** 535–547. MR2002227

[20] KIM, Y. and LEE, J. (2001). On posterior consistency of survival models. *Ann. Statist.* **29** 666–686. MR1865336

[21] KINGMAN, J. F. C. (1967). Completely random measures. *Pacific J. Math.* **21** 59–78. MR0210185

[22] KINGMAN, J. F. C. (1993). *Poisson Processes. Oxford Studies in Probability* **3**. Clarendon Press, New York. MR1207584

[23] LIJOI, A., PRÜNSTER, I. and WALKER, S. G. (2008). Posterior analysis for some classes of nonparametric models. *J. Nonparametr. Stat.* **20** 447–457. MR2424252

[24] LO, A. Y. and WENG, C.-S. (1989). On a class of Bayesian nonparametric estimates. II. Hazard rate estimates. *Ann. Inst. Statist. Math.* **41** 227–245. MR1006487

[25] NIETO-BARAJAS, L. E. and WALKER, S. G. (2004). Bayesian nonparametric survival analysis via Lévy driven Markov processes. *Statist. Sinica* **14** 1127–1146. MR2126344

[26] NIETO-BARAJAS, L. E. and WALKER, S. G. (2005). A semi-parametric Bayesian analysis of survival data based on Lévy-driven processes. *Lifetime Data Anal.* **11** 529–543. MR2213503

[27] PECCATI, G. and PRÜNSTER, I. (2008). Linear and quadratic functionals of random hazard rates: An asymptotic analysis. *Ann. Appl. Probab.* **18** 1910–1943.

[28] PECCATI, G. and TAQQU, M. S. (2008). Central limit theorems for double Poisson integrals. *Bernoulli* **14** 791–821.

[29] PETERSON, A.V. (1977). Expressing the Kaplan–Meier estimator as a function of empirical subsurvival functions. *J. Amer. Statist. Assoc.* **72** 854–858. MR0471165

[30] REGAZZINI, E., LIJOI, A. and PRÜNSTER, I. (2003). Distributional results for means of random measures with independent increments. *Ann. Statist.* **31** 560–585. MR1983542

[31] ROTA, G.-C. and WALLSTROM, C. (1997). Stochastic integrals: A combinatorial approach. *Ann. Probab.* **25** 1257–1283. MR1457619

[32] SATO, K. (1999). *Lévy Processes and Infinitely Divisible Distributions. Cambridge Studies in Advanced Mathematics* **68**. Cambridge Univ. Press, Cambridge. MR1739520

[33] SCHWARZ, L. (1965). On Bayes procedures. *Z. Wahrsch. Verw. Gebiete* **4** 10–26. MR0184378

[34] SURGAILIS, D. (1984). On multiple Poisson integrals and associated Markov semigroups. *Probab. Math. Statist.* **3** 217–239. MR0764148

[35]  WALKER, S. G. (2003). On sufficient conditions for Bayesian consistency. *Biometrika* **90** 482–488. MR1986664
[36]  WALKER, S. G. (2004). New approaches to Bayesian consistency. *Ann. Statist.* **32** 2028–2043. MR2102501
[37]  WU, Y. and GHOSAL, S. (2008). Kullback Leibler property of kernel mixture priors in Bayesian density estimation. *Electron. J. Stat.* **2** 298–331. MR2399197

P. DE BLASI
I. PRÜNSTER
DIPARTIMENTO DI STATISTICA
  E MATEMATICA APPLICATA
UNIVERSITÀ DEGLI STUDI DI TORINO
CORSO UNIONE SOVIETICA 218/BIS
10134 TORINO
ITALY
E-MAIL: pierpaolo.deblasi@unito.it
          igor@econ.unito.it

G. PECCATI
EQUIPE MODAL'X
UNIVERSITÉ PARIS OUEST NANTERRE LA DÉFENSE
200, AVENUE DE LA RÉPUBLIQUE
92000 NANTERRE
AND
LABORATOIRE DE STATISTIQUE
THÉORIQUE ET APPLIQUÉE
UNIVERSITÉ PARIS VI
FRANCE
E-MAIL: giovanni.peccati@gmail.com